\documentclass[11pt]{amsart}
\usepackage{amssymb,amsmath,amsthm}
\usepackage{epsfig}
\newtheorem{theorem}{Theorem}[section]
\newtheorem{lemma}[theorem]{Lemma}
\newtheorem{remark}[theorem]{Remark}
\newtheorem{definition}[theorem]{Definition}
\newtheorem{conjecture}[theorem]{Conjecture}
\newtheorem{proposition}[theorem]{Proposition}
\newtheorem{corollary}[theorem]{Corollary}

\newtheorem{example}[theorem]{Example}

\newcommand\A{{\mathcal A}}
\newcommand\B{{\mathcal B}}
\newcommand\M{{\mathbb M}}
\newcommand\Z{{\mathbb{Z}}}

\title{On the first group of the chromatic cohomology of graphs}
\author{Milena D. Pabiniak, J\'ozef H. Przytycki, Radmila Sazdanovi\'c}
\begin{document}
\maketitle 
\ \\
\ \\
\begin{quotation}
ABSTRACT. \baselineskip=10pt

The algebra of truncated polynomials $\A_m = \Z[x]/(x^m)$ plays an
important role in the theory of Khovanov and Khovanov-Rozansky
homology of links. 
 We have
demonstrated that Hochschild homology is closely related to
Khovanov homology via comultiplication free graph cohomology. 
It is not difficult to compute Hochschild
homology of $\A_m$ and the only torsion, equal to $\Z_m$, appears in
gradings $(i,\frac{m(i+1)}{2})$ for any positive odd $i$.
We
analyze here the grading of graph cohomology which is producing torsion
for a polygon. We find completely the cohomology
$H^{1,v-1}_{\A_2}(G)$ and $H^{1,2v-3}_{\A_3}(G)$. The group 
$H^{1,v-1}_{\A_2}(G)$ is closely related to the standard graph cohomology, 
except that the boundary of an edge is the sum of endpoints instead of 
the difference. 
The result about $H^{1,v-1}_{\A_2}(G)$ gives as a corollary a fact
about Khovanov homology of alternating and $+$ or $-$ adequate link
diagrams.
The group
$H^{1,2v-3}_{\A_3}(G)$ can be computed from the homology of a
cell complex, $X_{\Delta,4}(G)$, built from the graph $G$.
In particular, we prove that $\A_3$ cohomology can have any torsion.
We give a simple and complete characterization of those graphs which have
torsion in cohomology $H^{1,2v-3}_{\A_3}(G)$ (e.g. loopless graphs 
which have a 3-cycle).
We also construct graphs which have the same (di)chromatic 
polynomial but different $H^{1,2v-3}_{\A_3}(G)$.
Finally, we give examples of calculations of width of $H^{1,*}_{\A_3}(G)$ 
and of cohomology $H^{1,(m-1)(v-2)+1}_{\A_m}(G)$ for $m>3$.
\\
\end{quotation}
\ \\
\newpage

\tableofcontents
\section{Introduction}
In 1945, G.~Hochschild\footnote{According to \cite{Mac}, p.237:
"...From this result [that if a space has all homotopy groups, but
the first, $G$, trivial then homology of the space does not depend on
the choice of the space representative so is a property of the
group $G$ alone] developed the surprising idea that
cohomology, originally studied just for spaces, could also apply
to algebraic objects such as groups and rings. Given his
topological background and enthusiasm, Eilenberg was perhaps the
first person to see this clearly. He was in active touch with
Gerhard Hochschild, who was then a student of Chevally at
Princeton. Eilenberg suggested that there ought to be a cohomology
(and a homology) for algebras. This turned out to be the case, and
the complex used to describe the cohomology of groups (i.e., the
bar resolution) was adapted to define the Hochschild cohomology of
algebras. Eilenberg soon saw other possibilities for homology, and
he and Henri Cartan wrote the book Homological algebra, which
attracted lively interest among algebraists such as Kaplansky. A
leading feature was the general notion of a resolution, say of a
module $M$; such a resolution was an exact sequent of free modules
$F_i$ $... \to F_2 \to F_1 \to F_0 \to M\ $. Earlier work by
Hilbert on syzygies suggested this idea; an essential feature was
a theorem comparing two resolutions used to prove that the
cohomology they give is (up to isomorphism) independent of the
choice of the resolution".} defined a homology of any ring
\cite{Hoch}. The algebra of truncated polynomials $\A_m =
\Z[x]/(x^m)$  was one of the first to have its Hochschild homology
computed. This algebra plays an important role in the theory of
Khovanov and Khovanov-Rozansky homology of links. The only torsion
in the Hochschild homology $HH_{i,j}(\A_m)$ is equal to $\Z_m$ and
it appears in gradings $(i,\frac{m(i+1)}{2})$ for any odd $i$ (see
Section 2.1). We have demonstrated in \cite{Pr-2} that Hochschild
homology is closely related to Khovanov homology via
 comultiplication free graph cohomology described in
\cite{H-R-1,H-R-2}. In this paper we analyze graph cohomology in
gradings which have torsion for a polygon (as we know from the
relation with Hochschild homology). More precisely we are
interested in the group
$H^{1,(m-1)(v-1)-\frac{(m-2)(n-1)}{2}}_{\A_m}(G)$ for some odd $n$;
$v=v(G)$ denotes the number of vertices of $G$ (Section 2).
In the case of $n=3$ the corresponding cochain
complex has length two, \ $ C^{0,(m-1)(v-2)+1}
\stackrel{\partial^0}{\to}
  C^{1,(m-1)(v-2)+1} \to 0$,
and we find cohomology completely for $m=2$ and $m=3$, relating graph
cohomology to a cohomology of a cell complex built from the graph
$G$ (Theorem 4.1). 
For the case of $m=2$, which corresponds to classical Khovanov homology,
we obtain an interesting corollary about Khovanov homology of
alternating and adequate links (related to results of
\cite{Shu,A-P}). For the case of $m=3$ we show that the related cohomology
$H^{1,2v-3}_{\A_3}(G)$ can have arbitrary torsion and we characterize
those graphs which have torsion in $H^{1,2v-3}_{\A_3}(G)$. We illustrate
our results by several corollaries and by a few examples.
In the sixth section we construct graphs which have the same 
chromatic (even dichromatic and Tutte)
polynomial but different $H^{1,2v-3}_{\A_3}(G)$. In relation to these
examples we also analyze $A_3$ chromatic graph cohomology for one vertex
and two vertices product. We prove, in particular, that
$H^{1,2v-3}_{\A_3}(G)$ is a 2-isomorphism (matroid) graph invariant.
Finally, in the eight section, 
we report several computations  concerning width of
$tor(H^{1,*}_{\A_3}(G))$. We venture also into 
$H^{1,(m-1)(v-2)+1}_{\A_m}(G))$
calculations for $m>3$ and state several questions and conjectures.

\section{Basic facts about Hochschild homology and chromatic graph cohomology}

We recall in this section definitions of Hochschild homology and
chromatic graph cohomology and the relation among them \cite{Pr-2}.
 We follow \cite{Lo,H-R-2,H-P-R,Pr-2} in our exposition.

\subsection{Hochschild homology}

Let $k$ be a commutative ring and $\A$ a $k$-algebra (not
necessarily commutative).
Let $\M$ be a bimodule over $\A$, that is, a $k$-module on which
$\A$ operates linearly on the left and on the right in such a way
that $(am)a'=a(ma')$ for $a,a'\in \A$ and $m \in \M$. The actions
of $\A$ and $k$ are always compatible (e.g. $m(\lambda a)=
(m\lambda)a=\lambda (ma)$). When $\A$ has a unit element $1$ we
always assume that $1m=m1=m$ for all $m \in \M$. Under this unital
hypothesis, the bimodule $\M$ is equivalent to a right $\A\otimes
\A^{op}$-module via $m(a'\otimes a)=ama'$. Here $\A^{op}$ denotes
the opposite algebra of $\A$, that is, $\A$ and $\A^{op}$ are the
same as sets but the product $a\cdot b$ in $\A^{op}$ is the
product $ba$ in $\A$. The product map of $\A$ is usually denoted
$\mu: \A \otimes \A \to \A$,\ $\mu (a,b)=ab$.\\

In this paper we work only with unital algebras (the algebras of
truncated polynomials in most cases). We also assume, unless
otherwise stated, that $\A$ is a free $k$-module, however in most
cases, it suffices to assume that $\A$ is $k$-projective, or less
restrictively, that $\A$ is flat over $k$. \ Throughout the paper
the tensor product $\A \otimes \B$ denotes the tensor product over
$k$, that is, $\A \otimes_k \B$.

\begin{definition}[\cite{Hoch,Lo}]\label{Definition 2.1}\
The Hochschild chain complex $C_*(\A,\M)$ of the algebra $\A$
with coefficients in $\M$ is defined as:
$$ \ldots
\stackrel{b}{\to} \M\otimes \A^{\otimes n} \stackrel{b}{\to}
\M\otimes \A^{\otimes n-1} \stackrel{b}{\to} \ldots
\stackrel{b}{\to} \M\otimes \A \stackrel{b}{\to} \M$$ where
$C_n(\A,\M)= \M\otimes \A^{\otimes n}$ and the Hochschild boundary
is the $k$-linear map $b: \M\otimes \A^{\otimes n} \to \M\otimes
\A^{\otimes n-1}$ given by the formula
$b= \sum_{i=0}^n (-1)^id_i$, where the face maps $d_i$ are given by\\
$d_0(m,a_1,\ldots, a_n)= (ma_1,a_2,\ldots, a_n),$\\
$d_i(m,a_1,\ldots, a_n)= (m,a_1,\ldots , a_ia_{i+1},\ldots , a_n)$
for
$1\leq i \leq n-1$, \\
$d_n(m,a_1,\ldots, a_n)= (a_nm,a_1,\ldots ,a_{n-1})$.\\

In the case when $\M=\A$ the Hochschild complex is called the
{\it cyclic bar complex}.\\
By definition, the $n$th Hochschild homology group of the unital
$k$-algebra $\A$ with coefficients in the $\A$-bimodule $\M$ is
 the $n$th homology group of the Hochschild chain complex
denoted by $H_n(\A,\M)$. In the particular case $\M=\A$ we write
$C_*(\A)$ instead of $C_*(\A,\A)$ and $HH_*(\A)$ instead of
$H_*(\A,\A)$.
\end{definition}

The algebra $\A$ acts on $C_n(\A,\M)$ by $a\cdot
(m,a_1,...,a_n)=(am,a_1,...,a_n)$. If $\A$ is a commutative
algebra then the action commutes with boundary map $b$. Therefore,
$H_n(\A,\M)$ (in particular, $HH_*(\A)$) is an $\A$-module.

If $\A$ is a graded algebra, $\M$ a coherently graded
$\A$-bimodule, and the boundary maps are grading preserving, then
the Hochschild chain complex is a bigraded chain complex with ($b:
C_{i,j}(\A,\mathbb M) \to C_{i-1,j}(\A,\mathbb M))$, and
$H_{**}(\A,\mathbb M)$ is a bigraded $k$-module. In the case of
abelian $\A$ and $\A$-symmetric $\mathbb M$ (i.e. $am=ma$),
$H_{**}(\A,\mathbb M)$ is bigraded $\A$-module. The main examples
coming from the knot theory are $\A_m=\Z[x]/(x^m)$ and $\mathbb M$
the ideal in $\A_m$ generated by $x^{m-1}$. In the case of
$\M=\A_m$ we have (see for example \cite{Lo}).
\begin{proposition}\label{Proposition 2.2}
\begin{displaymath}
HH_{i,j}(\A_m)=
\left\{
\begin{array}{lll}
\Z_m & for   & i\ odd, \ j=\frac{i+1}{2}m\\
\Z & for & i=j=0\ \  or\ \ i\geq 0 \ \ and \\
& & \lfloor \frac{i}{2} \rfloor m +1 \leq
j \leq \lfloor \frac{i}{2} \rfloor m +m-1 \\
0 & otherwise &
\end{array}
\right.
\end{displaymath}
Here, $\lfloor x \rfloor$ denotes the integer part of $x$.\\
In particular, for $i$ odd $HH_{i,*}(\A_m)$ is $\A_m$-module
isomorphic to \\
$\Z[x]/(x^m,mx^{m-1}) \{m\frac{i-1}{2} +1\}$, where
$\{k\}$ denotes the shift by $k$ in the grading.
\end{proposition}

\subsection{Chromatic graph cohomology}\label{Subsection 2.2}
Chromatic graph cohomology was introduced in \cite{H-R-1} as a
comultiplication free version of Khovanov cohomology of
alternating links, where alternating link diagrams are translated
to plane graphs (Tait graphs). Being free of topological
restrictions, chromatic graph cohomology was extended in
\cite{H-R-2} to any commutative algebra $\A$. We showed in
\cite{Pr-2} that $\A$ graph cohomology (that is, chromatic graph
cohomology with underlining algebra $\A$) can be interpreted as a
generalization of Hochschild homology from a polygon to any graph.
We have this interpretation only for a commutative $\A$. It seems
to be, that if one works with general graphs and not necessary
commutative algebras then these algebras should satisfy some
''multiface" properties. Very likely planar algebras or
operads provide the proper framework.

\begin{definition}\label{Definition 2.3}
 For a given commutative $k$-algebra $\A$,
symmetric $\A$-module
$M$ and a graph $G$ with a base vertex $v_0$, we define
$M$-reduced $\A$ graph cochain complex and cohomology as follows (see
\cite{H-R-2,Pr-2} for details).
\begin{enumerate}
\item[(i)] The cochain  $k$-modules $C^i(G,v_0)$ are defined as follows:\\
$C^i(G,v_0) = \oplus_{|s|=i, s \subset E(G)} C^i_s(G,v_0) $. The
$k$-module $C^i_s(G,v_0) = M \otimes \A^{k(s)-1}$, where $k(s)$ is the
number of components of the graph $[G:s]$ which is the subgraph of
$G$ containing all vertices of $G$ and edges $s$.  We visualize
the product $M \otimes \A^{k(s)-1}$ as attachment of $M$ to a
component of $[G:s]$ containing $v_0$ and attachment of $\A$ to
any other component of $[G:s]$. \\ The cochain map $d^i :
C^i(G,v_0) \to C^{i+1}(G,v_0)$ is defined as follows: $d^i =
\sum_{e \notin s} (-1)^{t(s,e)} d_e^i$ where $d_e^i$ depends on
whether $e$ connects different components of $[G:s]$ or it connects
vertices in the same component of $[G:s]$. In the last case we
assume $d_e^i$ to be the identity map. If $e$ connects different
components of $[G:s]$ then either
\begin{enumerate}
\item[(m)] If $e$ connects the components of $[G:s]$ containing $v_0$
with another components, say the first one, then\\
$d_e^i(m,a_1,a_2,...a_{k(s)-1}) = (ma_1,a_2,...a_{k(s)-1})$,
\item[(a)] if $e$ connects two components not containing $v_0$,
say the first and the second, then\\
$d_e^i(m,a_1,a_2,...a_{k(s)-1}) = (m,a_1a_2,...a_{k(s)-1})$.
\end{enumerate}

\item[(ii)] We define $\M$-reduced cohomology denoted  by
$H^*_{\A,\mathbb M}(G,v_0)$ as the cohomology of the above cochain
complex. If we assume $\mathbb M = \A$ we obtain $\A$-cohomology
of graphs, $H^*_{\A}(G)$ (often called the chromatic graph cohomology as
it categorifies the chromatic polynomial of $G$ \cite{H-R-2}).
\end{enumerate}
\end{definition}
\begin{remark}\label{Remark 2.4}
 The boundary map $d_1: \M\otimes \A \to \M$ in 
Hochschild homology is the zero
map for a commutative algebra $\A$ and a symmetric module $\M$
($d_1(m,a)= ma-am=0$) and thus $H_0(\A,\M)=\M$.
Therefore, it is convenient to consider the variant of
chromatic graph cohomology, 
$\hat H^{*}_{\A,\M}(G,v_0)$, which has the zero map in place of Khovanov
comultiplication (see \cite{H-R-2,Pr-2}) as defined 
below\footnote{The modification allows us a concise formulation of Theorem 
2.5 even for a nonabelian $\A$.}:
\begin{enumerate}
\item[(i)] Consider the  cochain complex of a graph
$(\hat C^{*}_{\A,\M}(G,v_0); \hat d)$ obtained by modifying
$(C^{*}_{\A,\M}(G,v_0); d)$ as follows: $\hat C^i_{\A,\M}(G,v_0) =
C^{i}_{\A,\M}(G,v_0)$, $\hat d_e = d_e$ for a state $s$ such that $e$
is connecting different components of $[G:s]$ but if $e$ has
endpoints on the same component of $[G:s]$ we put $\hat d_e = 0$.
The cohomology of the cochain complex $(\hat C^{*}_{\A,\M}(G,v_0); \hat d)$
will be denoted by $\hat H^{*}_{\A,\M}(G,v_0)$. For $\M = \A$ we
write simply $\hat H^{*}_{\A}(G)$. This version of
chromatic graph cohomology was considered in \cite{H-R-2}.
\item[(ii)] Let $\ell(G)$ denote the girth of the graph $G$, that is,
the length of the shortest cycle in $G$. Then straight from the
definition of $H^{i}_{\A,\M}(G,v_0)$ and 
 $\hat H^{i}_{\A,\M}(G,v_0)$ we get:\\
$\hat H^{i}_{\A,\M}(G,v_0) = H^{i}_{\A,\M}(G,v_0)$ for $i < \ell(G)-1$, and
torsion part satisfies (for $k$ being a principal ideals domain)\\
$tor (\hat H^{i}_{\A,\M}(G,v_0)) = tor (H^{i}_{\A,\M}(G,v_0))$ for
$i = \ell(G)-1$.
\end{enumerate}
\end{remark}

\subsection{$\A$ graph cohomology of a polygon as Hochschild homology of $\A$}

We use the following result connecting Hochschild homology and
chromatic graph cohomology observed in \cite{Pr-2} (because Theorem 2.5
concerns only a polygon we can work with any, not necessary commutative,
algebra $\A$).
\begin{theorem}\label{Theorem 2.5} Let $\A$ be a unital
algebra which is a free $k$-module\footnote{We assume in
this paper that $\A$ is a free
$k$-module, but we could relax the condition to have $\A$ to be
projective or, more generally, flat over a commutative ring with
identity $k$; compare \cite{Lo}. We require $\A$ to be a unital
algebra in order to have an isomorphism $\M
\otimes_{\A^{\epsilon}} \A^{\otimes n+2} = \M \otimes \A^{\otimes
n}$; the isomorphism is given by $\M \otimes_{\A^{\epsilon}}
\A^{\otimes n+2} \ni(m,a_0,a_1,...,a_n,a_{n+1}) \to
(a_{n+1}ma_0,1,a_1,...,a_n,1)$ which we can write succinctly as
$(a_{n+1}ma_0,a_1,...,a_n) \in \M \otimes \A^{\otimes n}$. We
should stress that in $\M \otimes_{\A^{\epsilon}} \A^{\otimes
n+2}$ the tensor product is taken over $\A^{\epsilon}=\A\otimes
\A^{op}$ while in $\M \otimes \A^{\otimes n}$ the tensor product
is taken over $k$.}, $\mathbb M$ an $\A$-bimodule and $P_{n+1}$ --
the $(n+1)-gon$. Then for $0< i \leq n$ we have:
$$ \hat H^i_{\A,\mathbb M}(P_{n+1}) = H_{n-i}(\A,\mathbb M).$$
Furthermore, if ${\A}$ is a graded algebra and ${\mathbb M}$ a
coherently graded module then $ \hat H^{i,j}_{\A,\mathbb
M}(P_{n+1}) = H_{n-i,j}(\A,\mathbb M)$, for $0< i \leq n$ and
every $j$.
\end{theorem}
In particular,
 $ \hat H^{n,j}_{\A,\mathbb M}(P_{n+1}) =(sign?) H_{0,j}(\A,\mathbb M) = 
\M/(am-ma)$
\begin{corollary}\label{Corollary 2.6}
$ \hat H^{i,j}_{\A}(P_{n+1}) = HH_{n-i,j}(\A)$, for $0< i \leq n$
and every $j$. Furthermore, for a commutative $\A$,
$H^{i,j}_{\A}(P_{n+1}) = \hat H^{i,j}_{\A}(P_{n+1})$ , for $0< i
<n$ and $H^{n,*}_{\A}(P_{n+1}) = 0$, $\hat H^{n,*}_{\A}(P_{n+1})=
HH_{0,*}(\A)=\A$.\ For a general $\A$,  $\hat
H^{n,*}_{\A}(P_{n+1}) =HH_{0,*}(\A)=\A/(ab-ba)$.
\end{corollary}
We work in the paper with $\M=\A= \A_m$ but in the future analysis
for $M=(x^{m-1})$, $\A= \A_m$ will be given.

\subsection{Interesting gradings for $\A_m$-algebras}
The relation between Hochschild homology and graph cohomology 
of a polygon allowed us to find graph cohomology of a polygon for algebras
$\A_m$. In particular, the torsion of $H^{i,j}_{\A_{m}}(P_v)$ is 
supported by $(i,j)$ such that $v-i$ is even and $mi + 2j= mv$. We have:
\begin{corollary}[\cite{Pr-2}]\label{Corollary 2.7}
\begin{displaymath}
tor(H^{i,j}_{\A_m}(P_v)) =
\left\{
\begin{array}{ll}
H^{i,j}_{\A_m}(P_v) =\Z_m & for \  v-i\ even,    
\  0 < i \leq v-2, \ j=\frac{v-i}{2}m\\
0 & otherwise 
\end{array}
\right.
\end{displaymath}
\end{corollary}
The study of torsion in $H^{v-2,m}_{\A_m}(G)$ was initiated in \cite{H-P-R}.
In this paper we concentrate on the first cohomology $H^{1,j}_{\A_m}(G)$,
partially motivated by the fact that computing whole $H^{*,*}_{\A_m}(G)$ 
is NP-hard (so, up to famous conjecture, has exponential complexity) while 
computing $H^{1,j}_{\A_m}(G)$ for a fixed $m$ has polynomial complexity.

Corollary 2.7 (applied for $i=1$) suggest also that if a graph, $G$ has 
an odd cycle of length $n$ then the grading $j= (v-n)(m-1) + \frac{n-1}{2}m$ 
should be of considerable interest. We decided to work with the case 
$n=3$,  that is, to analyze $H^{1,(v-2)(m-1) +1}_{\A_m}(G)$.   

It is well known that $H^{1,j}_{\A_2}(G)$ is trivial for $j>v-1$, in fact 
the whole first cohomology is supported by $j=v-1,v-2$  
Also for $m=3$ the highest possible grade in
$H^{1,*}_{\A_3}$ is equal to $2v-3$,
 that is, $H^{1,j}_{\A_3}(G)=0$ for $j>2v-3$. This
gives another reason to concentrate on $j=2v-3$ grading.
We can ask whether for general $m$, the value $j=(v-2)(m-1) +1$ is the 
highest possible grading of nonzero $H^{1,j}_{\A_m}(G)$.
As a first step in this direction, sufficient for $m=3$, we 
prove the following proposition generalizing slightly
Corollary 13 (1c) of \cite{H-P-R}.
\begin{proposition}\label{Proposition 2.8}
Assume that  $m \geq 3$ and $0 < i < \ell (G)$, then
$H^{i,j}_{\A_m} = 0$ for $j \geq (m-1)(v-i)$
\end{proposition}
\begin{proof}
We use the following notation (following that of \cite{Vi-1}).
Any $s \in E(G)$ is called a state of $G$ and in the case of $\A=\A_m$,
an enhanced state $S$ is a state $s$ with every component of $[G:s]$
decorated by a weight $x^i$, $0 \leq i <m$.\
Notice first that for $i < \ell (G)$, and $j > (m-1)(v-i)$ we have
$C^{i,j}_{\A_3}(G) = 0$ and furthermore
$C^{i,(m-1)(v-i)}_{\A_3}(G)$ is freely generated by enhanced states
$S$, with underlining state $s$, and such that
each component of $[G:s]$ has weight $x^{m-1}$. For $i >0$ and $m>2$,
consider a component of $[G:s]$, say $X_e$ which has an edge $e$.
Consider the state $s'=s-e$ and the weights of components of $[G:s']$
in which one component of $X_e -e$ has weight $x$ and the other $x^{m-2}$.
Because $m-2 \geq 1$, therefore the image of this enhanced state is
the chosen generator $S$ of $C^{i,(m-1)(v-i)}_{\A_3}(G)$. Thus
$d_{i-1}$ is an epimorphism and $H^{1,(m-1)(v-i)}_{\A_3}=0$.
\end{proof}
We will discuss further improvements of Corollary 13 (1c) of \cite{H-P-R}
in the sequel paper (compare examples in Section 7). In particular we show
that $H^{1,j}_{\A_5}(G)=0$ for $j>4v-7$.

\subsection{From cohomology to homology}\label{Subsection 2.5}
In Sections 3 and 4 we are performing very concrete calculations and
we observe that it is much easier to work (visualize the chain complex)
in the homology case.
Homology and cohomology of a chain complex are related by the
universal coefficient theorem. We give here the simplified
version in the form we use (see for example \cite{Hat}).
\begin{proposition}\label{Proposition 2.9}
If the homology groups $H_n$ and $H_{n-1}$ of a chain complex
$C$ of free abelian groups are finitely generated then
$$H^n(C;\Z) = H_n(C;\Z)/tor(H_n(C;\Z)) \oplus tor(H_{n-1}(C;\Z)).$$
\end{proposition}
In particular, in the cases we consider mostly, we have: \\
$H^{0,(m-1)(v-2)+1}_{\A_m}(G) =
H_{0,(m-1)(v-2)+1}^{\A_m}(G)/tor(H_{0,(m-1)(v-2)+1}^{\A_m}(G))$, \\
$H^{1,(m-1)(v-2)+1}_{\A_m}(G) = H_{1,(m-1)(v-2)+1}^{\A_m}(G) \oplus
tor(H_{0,(m-1)(v-2)+1}^{\A_m}(G))$\\ 
and $H_{1,(m-1)(v-2)+1}^{\A_m}(G) = ker(C_{1,(m-1)(v-2)+1}^{\A_m}(G) \to 
C_{0,(m-1)(v-2)+1}^{\A_m}(G))$ is a free abelian group.

Because our (co)chain groups are free and finitely generated, we
use the same enhanced states to describe basis of chain,
$C_{i,j}^{\A_m}(G)$, and cochain, $C^{i,j}_{\A_m}(G)$, groups.
Matrices describing chain and cochain maps are transpose one to
another in these bases. In the cochain map we use multiplication
in algebra so in chain map we use dual comultiplication.
The concrete cases will be described in detail in next sections.


\section{The case of $\A_2=\Z[x]/(x^2)$ and Khovanov homology}\label{Section 3}

In this section we compute
$H^{1,v-1}_{\A_2}(G)$ for every graph, showing, in particular, that if $G$ is
connected then the torsion part of $H^{1,v-1}_{\A_2}(G)$ is either
trivial if $G$ is bipartite or otherwise (i.e. if $G$
has an odd cycle) it is equal to $\Z_2$.
In particular, the version of Shumakovitch's
conjecture for graphs holds for the height one (that is,
$H^{1,*}_{\A_2}(G)$ has no $\Z_4$ in its torsion part).
  The detailed analysis of $H^{2,v-2}_{\A_2}(G)$ will be given in the
sequel paper.

\begin{theorem}\label{A2-Jan1-2006 Theorem 3.1}
Let $G$ be a simple graph then
\begin{enumerate}
\item[(0)] $H^{0,v-1}_{\A_2}(G)= \Z^{p_0^{bi}}$, where $p_0^{bi}$ is
the number of bipartite components of $G$.
\item[(1)] $H^{1,v-1}_{\A_2}(G)= \Z^{p_1- (p_0 -p_0^{bi})} \oplus
\Z_2^{p_0 -p_0^{bi}}$, where $p_0$ is the number of components of
$G$ and $p_1= rank (H_1(G,\Z))= |E|-v+p_0$ is the cyclomatic number
of $G$.
\end{enumerate}
\end{theorem}
\begin{corollary}\label{Corollary 3.2}
If $G$ is a connected simple graph then
\begin{displaymath}
H^{1,v-1}_{\A_2}(G)=
\left\{
\begin{array}{ll}
\Z^{p_1} & if \ $G$ \ is\ a\ bipartite \ graph\\
\Z^{p_1 -1} \oplus \Z_2 &  if\ $G$\ has\ an\ odd\ cycle
\end{array}
\right.
\end{displaymath}

\end{corollary}
\begin{proof}
As mentioned in Section 2, it is easier to visualize the
chain complex in the homology case. That is, consider the
 interesting for us part of the chromatic graph chain complex
$$0 \leftarrow C^{\A_2}_{0,v-1}(G) \stackrel{d_1}{\leftarrow}
C^{\A_2}_{1,v-1}(G) \leftarrow 0.$$
We have $C^{\A_2}_{0,v-1}(G) = \Z^v$ and enhanced states forming
the basis of $C^{\A_2}_{0,v-1}(G)$ can be identified with vertices
of $G$, namely to an enhanced state in which every vertex but $v_i$ has 
attached the weight $x$ and $v_i$ has weight $1$ we associate the vertex 
$v_i$. We also have $C^{\A_2}_{1,v-1}(G) = \Z^E$ and enhanced states
forming the basis of $C^{\A_2}_{1,v-1}(G)$ can be identified with
edges of $G$ (we write $E$ for cardinality $|E|$ as long 
as the meaning is clear). The matrix describing the map $d_1$ is the incidence
matrix of the (unoriented) graph $G$ (see for example  \cite{Big}).
That is the image of an edge is equal to the sum of its endpoints.
Therefore, for a bipartite graph $H_{i,v-1}^{\A_2}(G)= H_i(G,\Z)$.
However, any odd cycle identifies a vertex on the cycle with its
opposite. It easily leads for a connected simple graph to
\begin{displaymath}
H_{0,v-1}^{\A_2}(G)=
\left\{
\begin{array}{ll}
\Z & if \ $G$ \ is\ a\ bipartite \ graph\\
\Z_2 &  if\ $G$\ has\ an\ odd\ cycle
\end{array}
\right.
\end{displaymath}
and for any simple graph to
$$H_{0,v-1}^{\A_2}(G)= \Z^{p_0^{bi}} \oplus \Z_2^{p_0 - p_0^{bi}}.$$
Expressing Euler characteristic in two ways we obtain:\\ 
 $rank(H_{0,v-1}^{\A_2}(G)) - rank(H_{1,v-1}^{\A_2}(G))=
rank(C_{0,v-1}^{\A_2}(G)) - rank(C_{1,v-1}^{\A_2}(G))= v-E$.\\
 Therefore,
for a connected simple graph $G$:
\begin{displaymath}
H_{1,v-1}^{\A_2}(G)=
\left\{
\begin{array}{ll}
\Z^{E-v+1}=\Z^{p_1} & if \ $G$ \ is\ a\ bipartite \ graph\\
\Z^{E-v}=\Z^{p_1-1} &  if\ $G$\ has\ an\ odd\ cycle
\end{array}
\right.
\end{displaymath}

and for any simple graph
$$H_{1,v-1}^{\A_2}(G)= \Z^{E-v + p_0^{bi}} =\Z^{p_1 - (p_0-p_0^{bi})} .$$
Corollary 3.2 and Theorem 3.1 follow almost immediately from the above 
results and Proposition 2.9.

\end{proof}

\subsection{From Kauffman states on link diagrams to graphs and
surfaces}\label{Subsection 3.1}

A Kauffman state $s$ of a link diagram $D$ is a function
from the set of crossings of $D$ to
the set $\{+1,-1\}$. Equivalently, to each crossing of $D$
we assign a marker
according to the following convention:\\
\ \\
\centerline{\psfig{figure=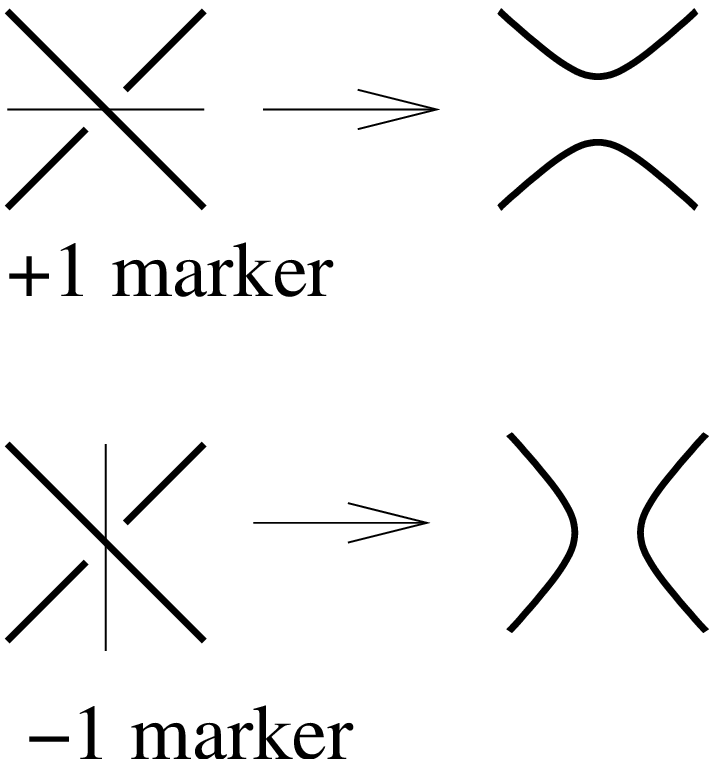,height=4.5cm}}
\begin{center}
Fig. 3.1; markers and associated smoothings
\end{center}

By $D_s$ we denote the system of circles in the diagram
obtained by smoothing all
crossings of $D$ according to the markers of the state $s$, Fig. 3.1.\\
By $|s|$ we denote the number of components of $D_s$.
The positive state $s_+$ (respectively the negative
state $s_-$) is the state with all positive markers
(resp. all negative markers).

For every Kauffman state $s$ of a link diagram $D$ we construct a
planar graph $G_s(D)$.
The graphs corresponding to states $s_+$ and $s_-$ are of particular
interest. If $D$ is an alternating diagram then $G_{s_+}(D)$ and
$G_{s_-}(D)$ are the plane graphs first constructed by Tait from
checkerboard coloring of regions of $R^2 - D$.

\begin{definition}\label{Definition 3.3}\ \
\begin{enumerate}
\item[(i)]
Let $D$ be a diagram of a link and $s$ its Kauffman state. We form
a graph, $G_s(D)$, associated to $D$ and $s$ as follows.
 Vertices of $G_s(D)$ correspond to circles of $D_s$.
Edges of $G_s(D)$ are in bijection with crossings of $D$ and an
edge connects given vertices if the corresponding crossing
connects circles of $D_s$ corresponding to the
vertices\footnote{If $S$ is an enhanced Kauffman state of $D$
then, in a similar manner, we associate to $D$ and $S$ the graph
$G_S(D)$ with signed vertices. Furthermore, we can additionally
equip $G_S(D)$ with a cyclic ordering of edges at every vertex
following the ordering of crossings at any circle of $D_s$. The
sign of each edge is the label of the corresponding crossing. In
short, we can assume that $G_S(D)$ is a ribbon (or framed) graph,
 and that with every state we associate a surface $F_s(G)$
whose core is the graph $G_s(D)$. $F_s(G)$ is naturally embedded
in $R^3$ with $\partial F_s(G)= D$. For $s=\vec{s}$, that is, $D$ is
oriented and markers of $\vec{s}$ agree with orientation of $D$,
$G_s(D)$ is the Seifert graph of $D$ and $F_s(G)$ is the Seifert
surface of $D$ obtained by Seifert construction. We do not use
this additional data in this paper but it may be of great use in
analysis of Khovanov homology (compare \cite{Pr-1}).}.
\item[(ii)]
In the language of associated graphs we can state the definition
of adequate diagrams as follows:\ the diagram $D$ is $+$-adequate
(resp. $-$-adequate) if the graph $G_{s_+}(D)$ (resp.
$G_{s_-}(D)$) has no loops.
\end{enumerate}
\end{definition}
In this language we can recall the result about torsion in
Khovanov homology \cite{A-P}; Theorem 2.2, which we generalize.
\begin{theorem}\cite{A-P}\label{Theorem 3.4} \ \\
Consider a link diagram $D$ of $N$ crossings. Then
\begin{enumerate}
\item[(+)]
If $D$ is $+$-adequate and $G_{s_+}(D)$ has a cycle of odd length,
then  the Khovanov homology has $\Z_2$ torsion.
More precisely, \\
$H_{N-2,N+2|s_+|-4}(D)$ has $\Z_2$ torsion,
\item[(-)] If $D$ is
$-$-adequate and $G_{s_-}(D)$ has a cycle
of odd length, then \\
$H_{-N,-N-2|s_-|+4}(D)$ has $\Z_2$ torsion.
\end{enumerate}
\end{theorem}

In \cite{H-P-R} we proved the following relation between graph
cohomology and classical Khovanov homology of alternating links.
\begin{theorem}\label{Theorem 3.5}
Let $D$ be the diagram of an unoriented framed alternating link
and let $G$ be its Tait graph (i.g. $G= G_{s_+}(D)$).
Let $\ell$ denote the girth of $G$, that is, the length of the
shortest cycle in $G$. For all $i<\ell-1$, we have

\[
H^{i,j}_{\A_2}(G) \cong H_{a,b}(D)
\]
with $\left
\{\begin{array}{l} a=E(G)-2i,\\
b=E(G)-2V(G)+4j,
\end{array}
\right.$

where $H_{a,b}(D)$ are the Khovanov homology groups of the
unoriented framed link defined by $D$, as explained in \cite{Vi-1}.

Furthermore, $tor(H^{i,j}_{\A_2}(G)) = tor(H_{a,b}(D))$ for
$i=\ell-1$.
\end{theorem}

We generalize Theorem 3.5 from an alternating diagram to any diagram
by using the graph $G_{s_+}(D)$ which for alternating diagrams is
a Tait graph. The proof follows exactly the same
line as that of Theorem 3.5 given in \cite{H-P-R}. Below we use notation
from Theorem 3.5.

\begin{theorem}\label{Theorem 3.6}
Let $D$ be the diagram of an unoriented framed link
and $G= G_{s_+}(D)$ its associated graph. Then:
\begin{enumerate}
\item[(i)]
For all $i<\ell-1$, we have
\[
H^{i,j}_{\A_2}(G) \cong H_{a,b}(D)
\]

\item[(ii)]
For $i=\ell-1$ we have  $tor(H^{i,j}_{\A_2}(G)) = tor(H_{a,b}(D))$.
\end{enumerate}
\end{theorem}

If girth $\ell(G_{s_+}(D))>2$ we say that $D$ is strongly $+$-adequate.
From the main result of this section (Theorem 3.1) and Theorem 3.6
 we get the following generalization of Theorem 3.4:
\begin{corollary}\label{Corollary 3.7}
\begin{enumerate}
\item[(i)] Assume that $D$ is a $+$-adequate diagram,
then \\
 $tor(H_{N-2,N+2|s_+|-4}(D)) = \Z_2^{p_0(G_{s_+}(D)-p_0^{bi}(G_{s_+}(D))}$.
\item[(ii)] Assume that $D$ is a strongly $+$-adequate diagram.
Then $$H_{N-2,N+2|s_+|-4}(D)) = \Z_2^{p_0-p_0^{bi}} \oplus
\Z^{p_1-(p_0-p_0^{bi})}.$$
\item[(iii)] Assume that $D$ is a $-$-adequate diagram,
then \\
$tor(H_{-N+2,-N-2|s_-|+4}(D)) = \Z_2^{p_0(G_{s_-}(D)-p_0^{bi}(G_{s_-}(D))}$.
\item[(iv)] Assume that $D$ is a strongly $-$-adequate diagram.
Then $$H_{-N+2,-N-2|s_-|+4}(D)) = \Z_2^{p_0-p_0^{bi}} \oplus
\Z^{p_1-(p_0-p_0^{bi})}.$$
\end{enumerate}
\end{corollary}
We can associate to any Kauffman state $s$ not only the
graph $G_{s}(D)$ but also a surface, $F_{s}(D)$, such that $G_{s}(D)$
is the spine of $F_{s}(D)$ (we generalize in such a way
Tait's black and white surfaces of checkerboard coloring and the Seifert
surfaces of an oriented diagram; compare Footnote 4.).
We can rephrase Corollary 3.7(i) to say
 that the torsion of $H_{N-2,N+2|s_+|-4}(D)$ is $\Z_2^n$ where
$n$ is the number of unoriented components of the surface $F_{s_+}(D)$.

We also speculate that there is a relation of $\A_m$ graph
homology to $sl(m)$ Khovanov-Rozansky \cite{K-R-1,K-R-2} homology and/or 
colored Jones homology of links \cite{Kh-2}.

To put Theorem 3.1 and Corollary 3.5 in perspective let us recall
that very little is known about torsion in Khovanov homology of
links. For a while it was thought that the only possible torsion
is 2-torsion (i.e. $\Z_2,\Z_4,\Z_8,...$). Then Bar-Natan announced that
torus knots can have odd torsion (e.g. homology of the torus knot
of type $(8,7)$ has $\Z_3,\Z_5$ and $\Z_7$ in its torsion
\cite{BN-3,BN-4}\footnote{For example the 22nd homology at degree 73 is
equal to $\Z_2 \oplus \Z_4 \oplus \Z_5 \oplus \Z_7$.}). For prime
alternating links the only torsion found so far is $\Z_2$ torsion.
A.~Shumakovitch proved that there is only $2$-torsion in Khovanov
homology of alternating links and he conjectures that
$\Z_4$-torsion is impossible\footnote{We use $p$-torsion to
describe elements of order $p^i$ for a prime $p$ while
$\Z_n$-torsion denotes an element of order $n$.}. Shumakovitch
proved also that every alternating link which is not disjoint or
connected sum of Hopf links or trivial links has $\Z_2$ torsion
\cite{Shu}. In \cite{A-P} we found explicity  $\Z_2$ torsion in
many adequate links, recovering in particular, the result of
Shumakovitch. In \cite{H-P-R} we proved the result from \cite{A-P}
in the $\A_2$ graph cohomology setting. In particular we proved that
a simple graph, which is not a forest has $\Z_2$ in
$H^{1,v-1}_{\A_2}(G)$ if $G$ has an odd cycle and it has $\Z_2$ in
$H^{2,v-2}_{\A_2}(G)$ if $G$ has and even cycle. In this section
we have computed completely $H^{1,v-1}_{\A_2}(G)$, showing, in
particular, that if $G$ is connected then the torsion part of
$H^{1,v-1}_{\A_2}(G)$ is trivial if $G$ is bipartite and it is
$\Z_2$ otherwise (i.e. $G$ has an odd cycle). In particular, the
version of Shumakovitch conjecture for graphs holds for height one
($H^{1,*}_{\A_2}(G)$ has no $\Z_4$ in its torsion part).

\section{Computation of $H^{1,2v-3}_{\A_3}(G)$}\label{Section 4}
The main result of this section describes the cohomology (and homology) 
at degree $2v-3$. The result is described using certain cell complex 
built from a graph $G$: 
$X_{\Delta,4}$ is the cell complex obtained from $G$ by adding 2-cells
along 4-cycles in $G$, identifying all vertices of $G$ and finally
adding 2-cells along expressions $2\vec e_3 - \vec e_2 - \vec e_1$
for any 3-cycle in $G$ -- two 2-cells added per every
3-cycle\footnote{To have uniquely defined cell complex we would have
 to add three 3-cells, but because then attachments would be
linearly dependent: $(2\vec e_3 - \vec e_2 - \vec e_1) +
(2\vec e_3 - \vec e_2 - \vec e_1) + (2\vec e_3 - \vec e_2 - \vec e_1) = 0$,
thus from the point of view of $H^1(X_{\Delta,4},\Z)$ or 
$H_1(X_{\Delta,4},\Z)$, it suffices to  
add two 2-cells. The different choice of two 2-cells for 
any 3-cycle of $G$ may however change the fundamental group 
of $X_{\Delta,4}$. This can be repaired by choosing relations of 
the form $\vec e_{i+1} - \vec e_{i} + \vec e_{i+1} - \vec e_{i-1}$, 
i=1,2,3, or in multiplicative notation 
$\vec e_{i}\vec e_{i-1}^{\,-1}\vec e_{i}\vec e_{i+1}^{\,-1}$. With this choice 
we have the identity $(\vec e_{3}\vec e_{2}^{\,-1}\vec e_{3}\vec e_{1}^{\,-1})
(\vec e_{1}\vec e_{3}^{\, -1}\vec e_{1}\vec e_{2}^{\,-1}) 
(\vec e_{2}\vec e_{1}^{\,-1}\vec e_{2}\vec e_{3}^{\,-1})=1$, thus 
$\pi_1(X_{\Delta,4})$ would not depend on the choice of two 2-cell 
attachments out of three possibilities. In the case of $G=P_3$, a triangle, 
we get $\pi_1((P_3)_{\Delta,4}) = \Z_3*\Z$.} (see Figure 4.9).
\begin{theorem}\label{Theorem 4.1}
 For an arbitrary simple 
 graph $G$ with $v$ vertices and
chromatic homology over algebra $\A_3$ the following is true:
 \begin{enumerate}
 \item [(0)] $H^{0,2v-3}_{\A_3}(G)$ is a free abelian group isomorphic to
  $H^1(X_{\Delta,4},\Z)
\oplus \Z^{t_0+\frac{d_2}{2}+d_{\geqslant 3}} $,\
 where $t_0$ is the number
of unoriented triplets of vertices not connected by any edge, $d_2$
is the number of ordered pairs of vertices of distance two and
$d_{\geqslant 3}$ is the number of ordered pairs of vertices of distance
at least three.
\item [(1)] $H^{1,2v-3}_{\A_3}(G) =
H^2(X_{\Delta,4},\Z) \oplus \Z^{t_2-\frac{d_2}{2}- sq(G)}$, 
 where\footnote{$t_2-\frac{d_2}{2}- sq(G)$ can be a negative number 
so formally it would be better to brake the formula into torsion and 
free part, that is: $tor H^{1,2v-3}_{\A_3}(G) = tor H^2(X_{\Delta,4},\Z)$ 
and $rank\, H^{1,2v-3}_{\A_3}(G) = rank\, H^2(X_{\Delta,4},\Z) + 
t_2-\frac{d_2}{2}- sq(G)$. Compare Remark 4.4.}
$t_2$ is the number of unoriented triplets of vertices connected
by exactly two edges (we call such a configuration, a joint), and 
$sq(G)$ denotes the number of squares (i.e. 4-cycles) in $G$.
\end{enumerate}
\end{theorem}
Below is the reformulation of Theorem 4.1 in the language of homology  which
will be used to  prove
our main result of this section.
\begin{theorem}\label{Theorem 4.2}
 For an arbitrary simple graph $G$ with $v$ vertices and
chromatic homology over algebra $\A_3$ the following is true:
 \begin{enumerate}
 \item [(0)]
$H_{0,2v-3}^{\A_3}(G)= H_1(X_{\Delta,4},\Z)
\oplus \Z^{t_0+\frac{d_2}{2}+d_{\geqslant 3}} $, \\
  \item [(1)] $H_{1,2v-3}^{\A_3}(G)$ is a free abelian group isomorphic to \\  
 $ H_2(X_{\Delta,4},\Z) \oplus \Z^{t_2-\frac{d_2}{2}- sq(G)}$.
 \end{enumerate}
\end{theorem}

 Since the proof of Theorem 4.2
requires a lot of technical details we will first
 give a brief outline containing main ideas and then the proof
 itself.
 As we have already mentioned in Section 2.5 
calculating homology
 instead of cohomology enables us to establish straightforward
 connections to the
homology of a cell complex corresponding to our graph.
  In order to get more information about  $H_{0,2v-3}^{\A_3}(G)$  and
$H_{1,2v-3}^{\A_3}(G)$
 we will calculate homology with coefficients in $\mathbb Z_3$ and
 $\mathbb Z [\frac{1}{3}]$ (i.e. the localization on the multiplicative
set generated by 3).
 In some special cases we will be able to distinguish $\Z_{3^i}$-
 torsion for an arbitrary $i$ by computing $3H_{0,2v-3}^{\A_3}(G)$.

 Recall that $ C_0(G)\cong \A_m^{\otimes v}$ and because $\A_m$ is a
 free abelian group with basis $1,x,x^2,...,x^{m-1}$, thus
 $C_0^{\A_m}$ has basis of
 $v$-tuples. In particular,\\
 $C_{0, (m-1)(v-2)+1}^{\A_m}(G)$ has
 basis  of $v$-tuples $(x^{a_1},x^{a_2},...,x^{a_v})$ with $0\leq
 a_i<m$ and $\sum_{i=1}^m a_i = (m-1)(v-2)+1$. 
Similarly, $C_1(G)\cong \oplus_E \A_m^{\otimes v-1}$ where $E$ will 
be used to denote set of edges in our graph
 $G$ and its cardinality, as long as it causes no confusion.
Because $G$ is a simple graph ($\ell(G)\geq 2$) therefore for $|s|=2$ the 
graph $[G:S]$ has $k(s)=v-2$ components and thus 
$C_{2,(m-1)(v-2)+1}^{\A_m} = 0$ and $H_1^{\A_m}(G)= \ker (d_1)$.
We work in this section with $m=3$ and grading $2v-3$. In this grading
our chain complex is as follows:
$$0 \stackrel{}{\leftarrow} C_{0,2v-3}(G) \stackrel{d_1}{\leftarrow}
C_{1,2v-3}(G) \stackrel{d_2}{\leftarrow} 0$$
The chain groups in this grading can be easily expressed
in terms of number of vertices $v$ and number of edges $E$ of our graph:\\

$C_{0,2v-3}^{\A_3}(G)$ has two essentially different types of
generators:
\begin{itemize}
   \item unordered triples $(v_i,v_j,v_k)=(x,x,x)$; vertices $v_i,v_j,v_k$ 
having weights equal to $x$; remaining vertices having weights  $x^2$
   \item ordered pairs $(v_i,v_j)=(1,x)$; weights of $v_i$ and $v_j$  
being equal to $ 1$ and $x$ respectively;
remaining vertices having weights  $x^2$
\end{itemize}
In particular  
$C_{0,2v-3}^{\A_3}(G) \cong  \mathbb{Z}^{(^v_3)+v(v-1)}.$

Similarly, $C_{1,2v-3}^{\A_3}$ has two types of generators
\begin{itemize}
   \item $(e_i)$ with weight $x$ and all isolated vertices have weights $x^2$
   \item $(e_i,v_j)$ where weight of the edge is $x^2$ and $v_j$ is the
   only vertex with weight $x$, the rest have weight $x^2$.
\end{itemize}
In particular, $ C_{1,2v-3}^{\A_3}(G)
\cong \mathbb{Z}^{E(v-1)}.$
 Recall that boundary map in cohomology uses multiplication in algebra,
 therefore differential in homology is its dual comultiplication.
 In particular, we consider all possible splits of components
 (in the case of $C_{1,2v-3}$ edge is mapped to its endpoints) and on
 the algebraic level this corresponds to all possible 
factorizations of a basic element $x^k$ representing weight of 
the particular component being split into $x^i$ and $x^{k-i}$.

Using the above notation for basis we describe decompositions of
$C_{0,2v-3}$ as a direct sum of free abelian groups:

$$C_{0,2v-3}(G)=C^{(t_0)}_{0,2v-3} \oplus
C_{0,2v-3}^{(t_1)}\oplus C_{0,2v-3}^{(t_2)}\oplus
C_{0,2v-3}^{(t_3)}\oplus C_{0,2v-3}^{(d_1)}\oplus
C_{0,2v-3}^{(d_2)}\oplus C_{0,2v-3}^{(d_3)}$$ where the summands are
defined as follows:
\begin{itemize}
 \item $C_{0,2v-3}^{(t_0)}$ is freely spanned by all triples $(v_i,v_j,v_k)$
 with no connections between vertices $v_i$,
$v_j$, $v_k$ in the graph $G$
 \item $C_{0,2v-3}^{(t_1)}$ is freely spanned by all triples $(v_i,v_j,v_k)$
 with exactly one of the edges between $v_i$,
$v_j$,$v_k$ present in the graph $G$
 \item $C_{0,2v-3}^{(t_2)}$ is freely spanned by all triples $(v_i,v_j,v_k)$
 with exactly $2$ of edges between vertices $v_i$,
$v_j$, $v_k$ belonging to the graph $G$
 \item $C_{0,2v-3}^{(t_3)}$ is freely spanned by all triples $(v_i,v_j,v_k)$
 with all $3$ edges between vertices $v_i$,$v_j$, $v_k$ present in the
 graph
 \item $C_{0,2v-3}^{(d_1)}$ is freely spanned by all ordered pairs $(v_i,v_j)$
with the distance $d(v_i,v_j)=1$
 \item $C_{0,2v-3}^{(d_2)}$ is freely spanned by all ordered pairs $(v_i,v_j)$
with the distance $d(v_i,v_j)=2$
 \item $C_{0,2v-3}^{(d_{\geq  3})}$ is freely spanned
by all ordered pairs $(v_i,v_j)$ with the distance $d(v_i,v_j)\geq 3$
 \end{itemize}
We use the convention that $t_0,t_1,t_2,t_3,d_1,d_2,d_{\geq 3}$
without parenthesis are the actual numbers and represent the ranks
of those groups.
\\ \ \\

Now we explain the direct sum decomposition of $C_{1,2v-3}$:
$$C_{1,2v-3}(G)=C^{(E)}_{1,2v-3} \oplus
C_{1,2v-3}^{(t_1)}\oplus C_{1,2v-3}^{(2t_2)}\oplus
C_{1,2v-3}^{(3t_3)}$$ where the following groups are subgroups of
$C_{1,2v-3}$, $e=\overline{v_i v_j}$ denotes and edge whose
endpoints are vertices $v_i$,$v_j$:
\begin{itemize}
  \item $C_{1,2v-3}^{(E)}$ is freely generated by states
  $(e)=(x)$.
  \item $C_{1,2v-3}^{(t_1)}$ is freely generated by
 states
  $(e,v_k)=(x^2,x)$ where $e$ is the only edge between vertices $v_i$,$v_j$,
  $v_k$ in the graph $G$.
  \item $C_{1,2v-3}^{(2t_2)}$ is freely generated
 by states $(e,v_k)=(x^2,x)$ where exactly two edges between
 $v_i$,$v_j$,$v_k$ are present in $G$.
   \item $C_{1,2v-3}^{(3t_3)}$ is freely generated
 by all states $(e,v_k)=(x^2,x)$ where all three edges between
 $v_i$,$v_j$,$v_k$ are present in $G$.
\end{itemize}

Now we can use these direct sum decomposition to analyze our chain
complex and extract its parts that give free part of homology.. \\

\textbf{Step 1.} 
According to the definition of the differential:
$$ \textrm{im}\,d_1 \subset \widehat{C^{(t_0)}_{0,2v-3}}\oplus
C_{0,2v-3}^{(t_1)}\oplus C_{0,2v-3}^{(t_2)}\oplus
C_{0,2v-3}^{(t_3)}\oplus $$ $$\oplus C_{0,2v-3}^{(d_1)}\oplus
C_{0,2v-3}^{(d_2)}\oplus C_{0,2v-3}^{(d_3)},$$
 where $\widehat C$ denotes the summand which is deleted from the sum.

Therefore $C_{0,2v-3}^{(t_0)}(G)$ contributes 
$\mathbb{Z}^{t_0}$ to homology $H_{0,2v-3}^{\A_3}(G)$.

\textbf{Step 2.}
 Now we will show that $d_1(C_{1,2v-3}^{(E)}(G))\subset
C_{0,2v-3}^{(d_1)}(G)$. As before, let $e=\overline{v_i v_j}$.
Then:
\begin{equation}
d_1(\overline{v_i v_j})= (1,x)+ (x,1)=\overrightarrow{e}+
\overleftarrow{e}
\end{equation}
In other words, the following equivalent relation holds
for every edge $e$ in $E$:
\begin{equation}
\overrightarrow{e}=-\overleftarrow{e}.
\end{equation}

We use this relation to reduce the size of chain complex by
$\frac{d_1}{2}=E$:
\begin{equation}
C_{0,2v-3}^{(d_1)}(G)/d_1(C_{1,2v-3}^{(E)}(G))=
C_{0,2v-3}^{(E)}(G)\cong \mathbb{Z}^E\cong
\Z^{2E}/(\overrightarrow{e}=-\overleftarrow{e}).
\end{equation}

\textbf{Step 3.} If exactly two vertices, say $v_i,v_j$, are
connected by an edge in the graph $G$, then we have the
isomorphism between $C_{0,2v-3}^{(t_1)}(G)$ and 
 $C_{1,2v-3}^{(t_1)}$. More precisely:
\begin{equation}
d_1(\overline{v_iv_j})=(x,x,x)+(1,x^2,x)+(x^2,1,x)
\end{equation}
This relation is used to present $(x,x,x)$ states by other
elements in a unique way.


\textbf{Step 4.} If exactly two edges are present in the graph G
we are analyzing:
$$ d_1 (C_{1,2v-3}^{(2t_2)})\subset C_{0,2v-3}^{(t_2)}(G) \oplus
C_{0,2v-3}^{(E)}(G) \oplus C_{0,2v-3}^{(d_2)}(G).$$ 
Without loss
of generality, let those two edges be denoted by
$\overline{v_iv_j}$,$\overline{v_jv_k}$; then we have: 
$$d_1(\overline{v_i v_j})=(x,x,x)+(1,x^2,x)+(x^2,1,x)$$
$$d_1(\overline{v_i v_k})=(x,x,x)+(1,x,x^2)+(x^2,x,1)$$ 
Thus we can eliminate states $(x,x,x)$ and be left with relation
\begin{equation}
(x,x^2,1)-(1,x^2,x)=(x^2,1,x)-(x,1,x^2)=
\overrightarrow{e_1}-\overleftarrow{e_2}
\end{equation}
This relation allows us to eliminate half of elements of 
$C_{0,2v-3}^{d_2}(G)$ and the other half  contributes $\Z^{d_2/2}$
to homology $H_{0,2v-3}^{\A_3}(G)$.
If there is another vertex $v_l$ from graph $G$ 
such that $\overline{v_i v_l},\overline{v_k v_l}\in E$ then, in a similar
way we have the relation:
\begin{equation}
(x,x^2,1)-(1,x^2,x)=(x^2,x,1)-(x,x^2,1)=
\overrightarrow{e_3}-\overleftarrow{e_4}
\end{equation}
Combining the above two relations using the common expression 
on the left side of each equation we get:
\begin{equation}
\overrightarrow{e_1}-\overleftarrow{e_2}+ 
\overrightarrow{e_3}-\overleftarrow{e_4}=0
\end{equation}
or equivalently:
\begin{equation}
\overrightarrow{e_1}+\overrightarrow{e_2}+\overrightarrow{e_3}+ 
\overrightarrow{e_4}=0
\end{equation}
This relation will be used later in the proof of Main Lemma.

\centerline{\psfig{figure=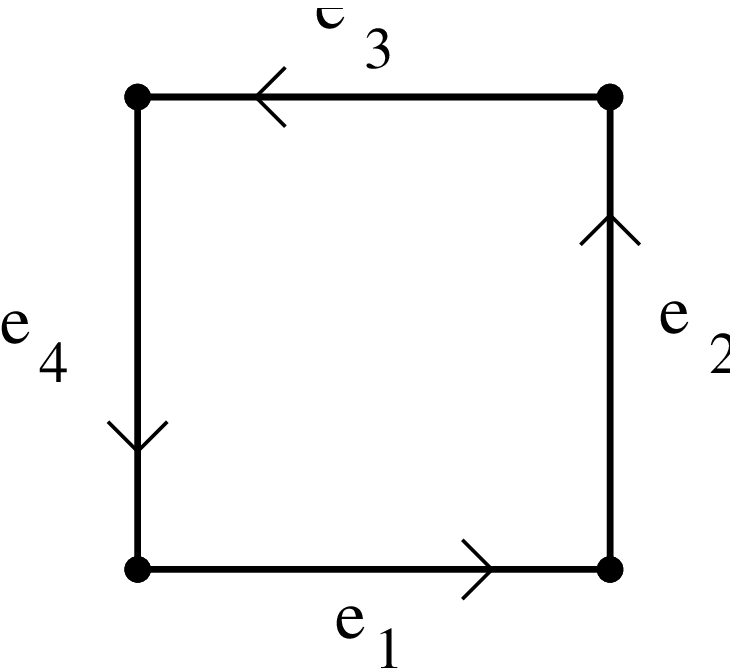,height=2.7cm}}
\centerline{ Figure 4.7}

 \textbf{Step 5.} Finally, if all three
edges are present in the graph G (forming a triangle in $G$)
then we have:

$$ d_1 (C_{1,2v-3}^{(3t_3)})\subset C_{0,2v-3}^{(t_3)}(G) \oplus
C_{0,2v-3}^{(E)}(G) $$ and as before we obtain following relations
(for simplicity, we denote the three edges of the triangle by
 $e_i, e_j, e_k$):\\
\begin{equation}
d_1(\overline{v_i
v_j})=(x,x,x)+(1,x^2,x)+(x^2,1,x)=\\
(x,x,x)+\overleftarrow{e_k}+\overrightarrow{e_j}
\end{equation}
\begin{equation}
d_1(\overline{v_jv_k})=\\(x,x,x)+(x,1,x^2)+(x,x^2,1)=\\
(x,x,x)+\overleftarrow{e_i}+\overrightarrow{e_k}
\end{equation}
\begin{equation}
d_1(\overline{v_i
v_k})=\\
(x,x,x)+(1,x,x^2)+(x^2,x,1)=(x,x,x)+\overrightarrow{e_i}+\overleftarrow{e_j}
\end{equation}
Therefore we can reduce the size of the chain group $C_{0,2v-3}$
 by eliminating $(x,x,x)$ states. After the elimination
 we are left with two relations:

\begin{equation}
\vec e_i- \vec e_j=\vec e_j-\vec e_k= \vec e_k- \vec e_i
\end{equation}
where edges of the triangle $v_i,v_j,v_k$ are coherently oriented
(Fig. 4.8).

\centerline{\psfig{figure=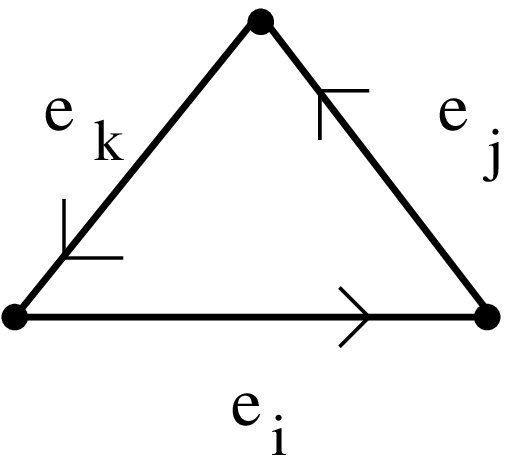,height=2.2cm}}
\centerline{ Figure 4.8: }


\textbf{Step 6.}
To formulate Main Lemma we need one more definition. Let $I_{\Delta,4}$
denote the subgroup of $\Z^{E}$ generated by two types of elements:\\
\begin{enumerate}
\item[(1)] Every square in our graph corresponds to one generator of $I_{\Delta,4}$:
$$u_4:=\overrightarrow{e_1}+ \overrightarrow{e_2}+\overrightarrow{e_3}+\overrightarrow{e_4}$$
\item[(2)] Every triangle from graph G contributes
 three generators to  $I_{\Delta,4}$:
 \\
$$ u_1:= (\overrightarrow{e_1}-\overrightarrow{e_2})- (\overrightarrow{e_2}-\overrightarrow{e_3})$$
$$ u_2:= (\overrightarrow{e_2}-\overrightarrow{e_3})- (\overrightarrow{e_3}-\overrightarrow{e_1})$$
$$ u_3:= (\overrightarrow{e_3}-\overrightarrow{e_1})- (\overrightarrow{e_1}-\overrightarrow{e_2})$$
Note that these generators are linearly dependent since \\
$u_1+u_2 + u_3 =0$.
\end{enumerate}

\centerline{\psfig{figure=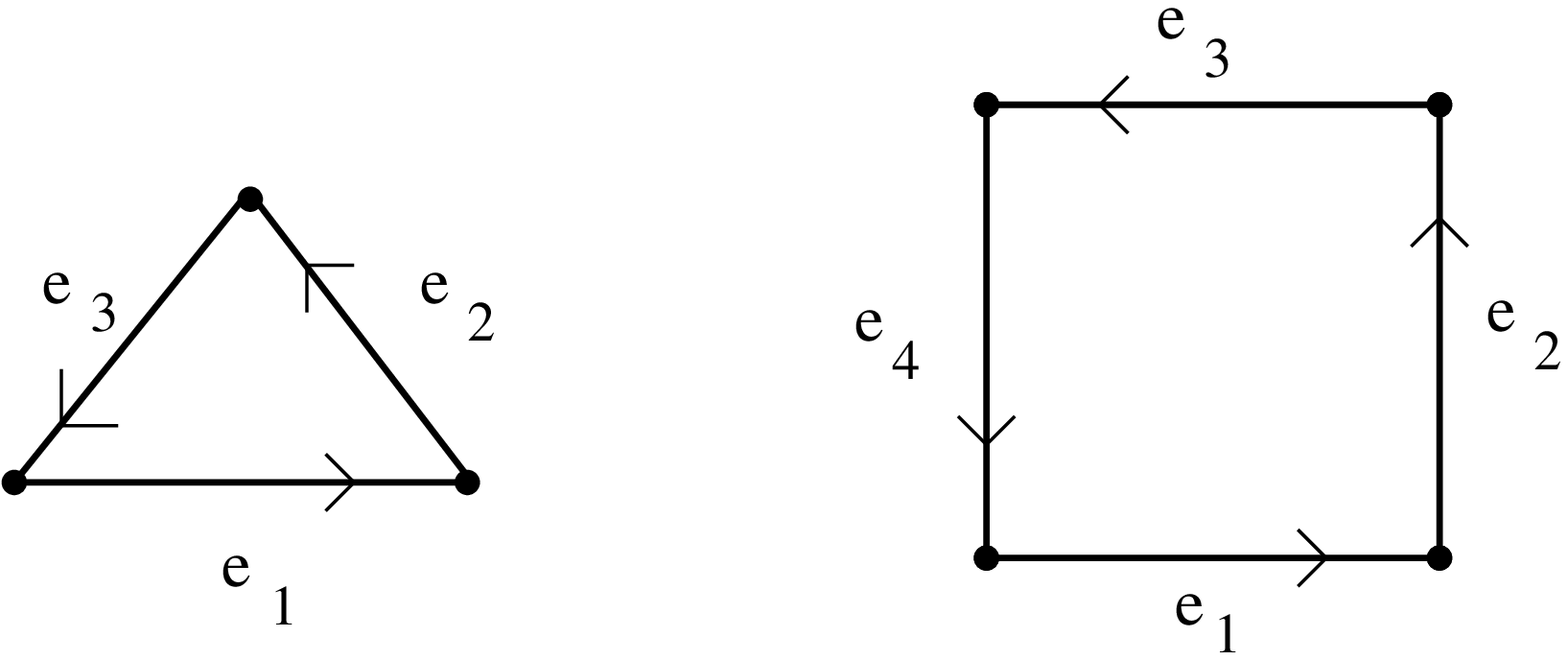,height=3.4cm}}
\centerline{ Figure 4.9:}
With this notation we can now summarize our analysis in
  Main Lemma.

\begin{lemma} [\textbf{Main Lemma}]\label{Lemma 4.3}

$$H_{0,2v-3}^{\A_3}(G)=\mathbb{Z}^E/I_{\Delta,4} \oplus
\Z^{t_0+\frac{d_2}{2}+d_{\geqslant 3}}.$$
\end{lemma}

\begin{proof} First,
$\Z^{t_0+\frac{d_2}{2}+d_{\geqslant 3}} $ is the free summand
of $H_{0,2v-3}^{\A_3}(G)$ by
Steps 1, 2, and by the fact that $im(d_1(C_{1,2v-3}(G)))$ belongs to
the direct sum of other then $\Z^{d_{\geqslant 3}}$ summands of 
$C_{0,2v-3}^{\A_3}(G)$.
 Recall that all generators of $C_{0,2v-3}$
are of the form $(x,x,x)$ or $(1,x)$ and in the case endpoints of 
an edge $e_i$ have weights $1$ and $x$ we write 
$(1,x)=\overrightarrow{e_i}$ or $(x,1)=\overleftarrow{e_i}$. 
In Steps 3, 4 and 5 we described how
to eliminate generators of type $(x,x,x)$.
In Step 2 we found  relation
$\overrightarrow{e}=-\overleftarrow{e}$. 
In Step 4 we also eliminated half of relations of the form $(1,x)$ 
when the distance between points is equal to two.
The other half contributes 
$\Z^{\frac{d_2}{2}}$ to $H_{0,2v-3}^{\A_3}(G)$. 
So finally we are left
with $E$ generators of type $\overrightarrow{e}$ spanning $\Z^E$.
Still, we have relations coming from squares and triangles (Steps 4 and 5).
These are exactly the relations generating $I_{\Delta,4}$.

Hence, $H_{0,2v-3}^{\A_3}(G)=\mathbb{Z}^E/I_{\Delta,4} \oplus
Z^{t_0+\frac{d_2}{2}+d_{\geqslant 3}}.$ 
\end{proof}
We are ready now to prove Theorems 4.1 and 4.2. 
\begin{proof}
We start from Theorem 4.2(0). Observe that $H_1(X_{\Delta,4};\Z)$ is 
equal to $\mathbb{Z}^E/I_{\Delta,4}$ because the cell complex 
$X_{\Delta,4}$ has as one skeleton the graph $G$ with all vertices 
identified so 1-cycles have $E$ as a basis. Furthermore, 2-cells 
of $X_{\Delta,4}$ were chosen in such a way that their boundaries
 generate the subgroup $I_{\Delta,4}$. Therefore Theorem 4.2(0) follows 
from Main Lemma (Lemma 4.3).  \ Theorem 4.2(1) follows from the fact that 
$H_{1,2v-3}^{\A_3}(G) = ker(d_1: C_{1,2v-3} \to C_{0,2v-3})$ so it is 
a free abelian group of the rank equal to 
$$rank\, C_{1,2v-3} - rank\, C_{0,2v-3} + rank\,H_{0,2v-3}^{\A_3}(G)$$
$$ =
E(v-1) - {v \choose 3} - v(v-1)+ rank\,H_{0,2v-3}^{\A_3}(G).$$
Furthermore
 $$rank\, H_{0,2v-3}^{\A_3}(G) =
rank\,H_1(X_{\Delta,4}(G)) + t_0 + \frac{d_2}{2}+d_{\geqslant 3}$$  and 
$$rank\,H_1(X_{\Delta,4}(G))=
 rank\,H_2(X_{\Delta,4}(G)) + rank\,H_0(X_{\Delta,4}(G)) 
+E -2t_3 -sq(G)-1$$ 
$$= rank\,H_2(X_{\Delta,4}(G)) + E -2t_3 -sq(G).$$
  Combining these together we get:
\\ $rank\, H_{1,2v-3}^{\A_3}(G) = $\\ 
$rank\,H_2(X_{\Delta,4}(G)) + 
E(v-1) - {v \choose 3} - v(v-1) + t_0 + \frac{d_2}{2}+d_{\geqslant 3}
+ E -2t_3 -sq(G)$\\
$= rank\,H_2(X_{\Delta,4}(G)) + t_2 - \frac{d_2}{2} - sq(G) + $\\
$(-t_2 + E(v-1) - {v \choose 3} - v(v-1) +t_0 + d_2 +d_{\geqslant 3} + E -2t_3)$ \\
$= rank\,H_2(X_{\Delta,4}(G)) + t_2 - \frac{d_2}{2} - sq(G)$,\\ as required. 
The last equality follows from the identities:
$$v(v-1)= d_1 +d_2 + d_{\geqslant 3} = 2E +d_2 + d_{\geqslant 3},$$  
$$E(v-1) = E + t_1 +2t_2 + 3t_3,$$  
$${v \choose 3} = t_0 +t_1 +t_2 +t_3,$$  and therefore 
$$-t_2 + E(v-1) - {v \choose 3} - v(v-1) +t_0 + d_2 +d_{\geqslant 3} + 
E -2t_3 =0.$$

Theorem 4.1 follows from Theorem 4.2 by applying Proposition 2.9.  
\end{proof}
 
\begin{remark}\label{Remark 4.4}
If $G$ has a 4-cycle with a diagonal (Figure 4.10) then 
the relation $0 = \vec{e_1} + \vec{e_2} + \vec{e_3} + \vec{e_4}$ 
yielded by this 4-cycle 
can be obtained from relations obtained from it follows from triangle
relations associated to triangles dividing the 
4-cycle. Namely, 
triangles give relations: $2\vec{e} = \vec{e_1} + \vec{e_2}$ and 
$-2\vec{e} = \vec{e_3} + \vec{e_4}$, whose sum is exactly the relation from 
the 4-cycle. Consequently, it may be useful to consider 
the cell complex $X_{\Delta,4'} \subset X_{\Delta,4}$ 
in which we glue $2$-cells along 4-cycles only in the case 
in which the 4-cycle has no diagonal. We observe that 
$H_1(X_{\Delta,4'},\Z) = H_1(X_{\Delta,4},\Z)$ and 
$H_2(X_{\Delta,4'},\Z) \oplus \Z^{sq'(G)} = H_1(X_{\Delta,4},\Z)$, 
where $sq'(G)$  denotes the number of $4$-cycles in $G$ which have a 
diagonal. In this notation Theorem 4.1(1) has the form:
$$H^{1,2v-3}_{\A_3}(G) =
H^2(X_{\Delta,4'},\Z) \oplus \Z^{t_2-\frac{d_2}{2}- sq(G) + sq'(G)}.$$
In Lemma 4.7 we will consider graphs in which every 
4-cycle has a diagonal, called square cordial (e.g. complete graphs 
or wheels). In this case $X_{\Delta,4'}= X_{\Delta}$  (no 2-cell is 
attached to a 4-cycle) which allows simpler formulation of 
Theorems 4.1 and 4.2.
\end{remark}
\centerline{\psfig{figure=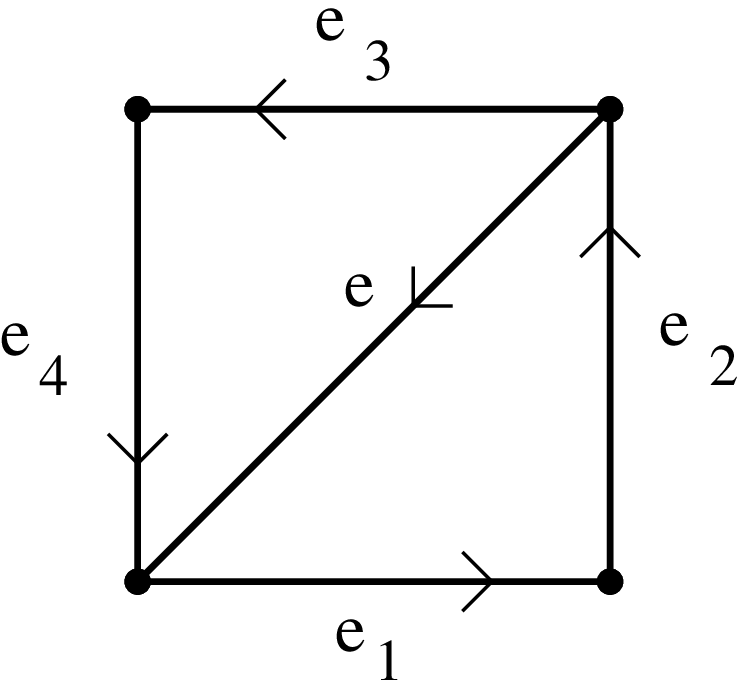,height=3.4cm}}
\centerline{ Figure 4.10:}

The case of square cordial graphs will be revisited in Lemma 4.7 and 
Corollary 4.8. 

Computing homology of
the cell complex 
$X_{\Delta,4}(G)$ is cumbersome, however 
we can recover   
a substantial part of it by considering three simpler complexes 
(compare Footnote 11):\\
(i)  $X_{(3),4}$ obtained from $G$ by identifying edges of every triangle 
in a coherent way (i.e. $\vec e_i = \vec e_j =\vec e_k$ for a triangle oriented 
as in Figure 4.8), then adding 2-cells along every 4-cycle in $G$ and 
finally identifying all vertices of $G$.\\
(ii) $\hat X_{3,4}$ obtained from $G$ by adding 2-cells 
along every 3- and 4-cycle of $G$.\\
(iii) 
$X_{3,4}$ obtained from $\hat X_{3,4}$ by identifying all vertices of $G$.\\
Observe that 
$H_2(X_{3,4}) = H_2(\hat X_{3,4}), 
H_1(X_{3,4}) = H_1(\hat X_{3,4}) \oplus \Z^{v-1}$ and \\
$H_1(X)\otimes R = H_1(X,R)$ according to 
the Universal Coefficient Theorem (e.g. \cite{Hat}), (compare also proof of Proposition 4.5).

\begin{proposition}\label{Proposition 4.5}\ 
\begin{enumerate}
\item[(a)] $H_1(X_{\Delta,4})\otimes \Z_3 = H_1(X_{\Delta,4},\Z_3) = 
 H_1(X_{3,4},\Z_3) =  H_1(X_{3,4})\otimes \Z_3$.
\item[(b)] 
$H_1(X_{\Delta,4})\otimes \Z[\frac{1}{3}] = 
H_1(X_{(3),4})\otimes \Z[\frac{1}{3}] = H_1(X_{(3),4},\Z[\frac{1}{3}])$.
\end{enumerate}
\end{proposition}
\begin{proof}
We know that  $H_1(X_{\Delta,4},\Z) =\Z^E/I_{\Delta,4}$
so we will consider the following short exact
sequence:
$$ C : 0 \rightarrow I_{\Delta , 4} \rightarrow \mathbb{Z}^E
\rightarrow \mathbb{Z}^E/I_{\Delta,4} \rightarrow 0 $$ From right
exactness of tensor product we get that the sequence:
$$I_{\Delta , 4} \otimes R \rightarrow \mathbb{Z}^E \otimes R
\rightarrow \mathbb{Z}^E/I_{\Delta,4}  \otimes R \rightarrow 0 $$ is
exact for every ring $R$. Hence
$$ (\mathbb{Z}^E/I_{\Delta,4})  \otimes R =(\mathbb{Z}^E\otimes R)/(I_{\Delta,4}  \otimes
R)= R^E/I^R_{\Delta,4}.$$ 
where  $I^R_{\Delta,4}$ is a submodule of $R^E$ constructed in the same way 
as $I_{\Delta,4}=I^{\Z}_{\Delta,4}$.

 In particular, 
$H_1(X_{\Delta,4})\otimes R = H_1(X_{\Delta,4},R)$. 
To prove Proposition 4.5 we consider cases 
 $R=\mathbb{Z}_3$ and
$R=\mathbb{Z}\bigl[\frac{1}{3}\bigr]$ separately.
\ \\
(a) Let $R=\mathbb{Z}_3$. Then $H_1(X_{\Delta,4},\Z_3) = 
  \Z_3^E/I^{\Z_3}_{\Delta,4}$
\\
To find $I^{\Z_3}_{\Delta,4}$ notice that in $\Z_3$ we have:
$$ e_1+e_3-2e_2=e_1+e_2+e_3 -3 e_2=e_1+e_2+e_3.$$
Similarly $e_2 + e_1 - 2e_2 =e_1+e_2+e_3 = e_3+e_2 -2e_1$.
 Thus $I^{\Z_3}_{\Delta , 4}$ is generated by 
expressions coming only from $3$ and $4$-cycles:
\begin{equation}
\overrightarrow{e_1}+\overrightarrow{e_2}+\overrightarrow{e_3}
\end{equation}
\begin{equation}
\overrightarrow{e_1}+\overrightarrow{e_2}+
\overrightarrow{e_3}+\overrightarrow{e_4}
\end{equation}

In this way we obtain that $\Z_3^E/I^{\Z_3}_{\Delta,4}$ is exactly 
the first homology of the cell complex $X_{3,4}(G)$
 with $\mathbb{Z}_3$ coefficients. That is, $\Z_3^E/I^{\Z_3}_{\Delta,4} = 
H_1(X_{3,4}(G),\Z_3) = H_1(X_{3,4}(G),\Z) \otimes \Z_3$.
\\
(b) Let $R=\Z[\frac{1}{3}]$. 
Then $H_1(X_{\Delta,4},\Z[\frac{1}{3}])=
\Z[\frac{1}{3}]^E/I^{\Z[\frac{1}{3}]}_{\Delta,4}$
\\
To find $I^{\Z[\frac{1}{3}]}_{\Delta,4}$ notice that 
because $3$ is invertible in  $\Z[\frac{1}{3}]$ 
therefore generators of $I^{\Z[\frac{1}{3}]}_{\Delta,4}$ coming 
from 3-cycles in $G$ are $e_i - \frac{e_1 + e_2 + e_3}{3}$, $i=1,2,3$
(they yield $e_1=e_2=e_3$ in the quotient).
Other generators are coming from 4-cycles and as before are of the form 
$\overrightarrow{e_1}+\overrightarrow{e_2}+
\overrightarrow{e_3}+\overrightarrow{e_4}$.\\
In this way we conclude that $\Z[\frac{1}{3}]^E/I^{\Z[\frac{1}{3}]}_{\Delta,4}$ 
is exactly the first homology of the cell complex $X_{(3),4}(G)$
 with $\Z[\frac{1}{3}]$ coefficients. That is, 
$\Z[\frac{1}{3}]^E/I^{\Z[\frac{1}{3}]}_{\Delta,4} =
H_1(X_{(3),4}(G),\Z[\frac{1}{3}]) = 
H_1(X_{(3),4}(G),\Z) \otimes \Z[\frac{1}{3}]$.

\end{proof}
In order to utilize and generalize Proposition 4.5 
recall that every finitely generated abelian group $H$ can be decomposed
uniquely:
$$H= \textrm{free}(H) \oplus \textrm{tor}(H)= \mathbb{Z}^a \oplus
\textrm{tor}_3(H) \oplus \textrm{tor}_{(3)}(H)$$ where
$\textrm{tor}_3(H)$ contains summands of the form
$\mathbb{Z}_{3^i}$, $i\in \{1,2,...\}$ and $\textrm{tor}_{(3)}(H)$
contains summands of the form $\mathbb{Z}_{n}$, where
$\textrm{gcd}(n,3)=1$.
\\
If we find $H_{(3)}\stackrel{def}{=}H \otimes \mathbb{Z}[\frac{1}{3}]$
we recover exactly $\textrm{free}(H) \oplus \textrm{tor}_{(3)}(H)$
as $$H \otimes \mathbb{Z}[\frac{1}{3}]=
(\mathbb{Z}[\frac{1}{3}])^a \oplus \textrm{tor}_{(3)}(H).$$
If $\textrm{tor}_3(H)=\mathbb{Z}^{a_1}_3 \oplus
\mathbb{Z}^{a_2}_{3^2} \oplus \ldots \oplus
\mathbb{Z}^{a_j}_{3^j}$ then 
$$H \otimes \mathbb{Z}_3=\mathbb{Z}_3^a \oplus
(\textrm{tor}_{(3)}(H) \otimes \mathbb{Z}_3)=\mathbb{Z}_3^a \oplus
\mathbb{Z}_3^{a_1+a_2+\ldots + a_j}.$$  Therefore $H_{(3)}$ and 
$H_3$ allows us to find 
$\Z^a$ and $\textrm{tor}_{(3)}(H)$ part
as well as the sum of exponents  $a_1+a_2+\ldots + a_j$ in
$\textrm{tor}_3(H)$ part.
Of course these is not sufficient to distinguish $\Z_{3^s}$ from
$\Z_3$. Sometimes, however, we can get enough information of 
$ 3H= \mathbb{Z}^a \oplus (\mathbb{Z}^{a_2}_3 \oplus
\mathbb{Z}^{a_3}_{3^2} \oplus \ldots \oplus
\mathbb{Z}^{a_j}_{3^{j-1}}) \oplus \textrm{tor}_{(3)}(H)$ 
to compute whole $H$. We illustrate it by first generalizing slightly 
Proposition 4.5 and then computing homology for 
an important class of graphs (square-cordial graphs),
including the complete graphs, $K_n$ and wheels, $W_n$ 
(that is, cones over ($n-1$)-gons, Figure 8.3), 
in which cases $3H$ has no 3-torsion.

\begin{proposition}\label{Proposition 4.6}
 Let $G$ be a simple graph and 
$$H_1(X_{\Delta,4}) = \mathbb{Z}^a \oplus
\textrm{tor}_3(H_1(X_{\Delta,4})) \oplus 
\textrm{tor}_{(3)}(H_1(X_{\Delta,4}))$$ 
where $\textrm{tor}_3(H_1({X(\Delta,4}))= \mathbb{Z}^{a_1}_3 \oplus
\mathbb{Z}^{a_2}_{3^2} \oplus \ldots \oplus
\mathbb{Z}^{a_j}_{3^j}$ then \\
(i) $ rank H_1(X_{(3),4},\Z[\frac{1}{3}])= rank H_1(X_{(3),4})= 
rank 3H_1(X_{\Delta,4}) = a$\\
(ii) $\textrm{tor}_{(3)}(H_1(X_{\Delta,4})) = 
\textrm{tor}_{(3)}(H_1(X_{(3),4},\Z[\frac{1}{3}]))= 
\textrm{tor}_{(3)}(H_1(X_{(3),4}))= \textrm{tor}_{(3)}3(H_1(X_{\Delta,4}))$, \\
(iii) $dim (H_1(X_{(3),4},\Z_3) = a +a_1 +a_2 + \ldots + a_j$\\
(iv) $\textrm{tor}_33H_1(X_{\Delta,4}) = \mathbb{Z}^{a_2}_3 \oplus
\mathbb{Z}^{a_3}_{3^2} \oplus \ldots \oplus
\mathbb{Z}^{a_j}_{3^{j-1}}$. Furthermore $\textrm{tor}_33H_1(X_{\Delta,4})$ 
is the quotient of $\textrm{tor}_3H_1(X_{(3),4})$.\\
In particular, if 
$H_1(X_{(3),4})$ has no 3-torsion then $H_1(X_{\Delta},4)$ has
no $Z_9$-torsion and $H_1(X_{\Delta,4})$ is fully determined by
$H_1(X_{(3),4})$ and $H_1(X_{3,4},\Z_3)$.
\end{proposition}
\begin{proof}
Parts (i),(ii), (iii) are just the reformulation of Proposition 4.5.
The first part of (iv) is obvious for any finitely generated group. 
The second part of (iv) follows from the fact that the equality 
$3\vec e_i=3 \vec e_j =3\vec e_k$ holds for any edges of coherently 
oriented triangle (Figure 4.8) and edges $3e$ generate $3H_1(X_{\Delta},4)$.
$3H_1(X_{\Delta,4})$ can have more relations than $H_1(X_{(3),4})$ thus 
we have the epimorphism from $H_1(X_{(3),4})$ to $3H_1(X_{\Delta,4})$ 
sending $e$ to $3e$ for any edge $e$. Because both groups have $\Z^a$ as 
a free part, therefore the 3-torsion part of $3H_1(X_{\Delta,4})$ 
is equal to a quotient of the 3-torsion part of $H_1(X_{(3),4})$. 
\end{proof}

A graph $G$ is called square-cordial if every $4$-cycle in
$G$ has a diagonal.
In order to be able to formulate our result about homology of square-cordial 
graphs, we define a new graph denoted by $G_{\Delta}$ and called the graph of
triangles of $G$.
The graph  $G_{\Delta}$ is obtained from $G$ as follows:
vertices of  $G_{\Delta}$ are in bijection with $3$-cycles of $G$.
Two vertices of $G_{\Delta}$ are connected by an edge if corresponding
triangles share an edge (we may assume that if the triangles share two edges
than vertices are connected by two edges). We say that the component 
of the graph $G_{\Delta}$ is coherent if relations $\vec e_i = 
\vec e_j = \vec e_k$ for every triangle in the component 
(compare Figure 4.8) never lead to relation of type 
$\overrightarrow{e} = \overleftarrow{e}$ (i.e. $2\vec e =0$).

\begin{lemma}\label{Lemma 4.7}
 Let $G$ be a simple square-cordial graph, then
$$H_1(X_{\Delta}) = \Z^{p_0(G_{\Delta})^{coh}} \oplus 
\Z_2^{p_0(G_{\Delta})-p_0(G_{\Delta})^{coh}} \oplus 
\Z_3^{v-1 + dim H_1(\hat X_3,\Z_3)- p_0(G_{\Delta})^{coh}},$$ where 
$p_0(G_{\Delta})^{coh}$ is the number of coherent components of $G_{\Delta}$.\\
In particular $H_1(X_{\Delta})$ has no $\Z_9$-torsion.
\end{lemma}
\begin{proof}
If a triangle graph $G_{\Delta}$ is connected, then clearly $$H_1(X_{(3)})= 
\begin{cases}  \Z   & if\ G_{\Delta} \ is\  coherent \\
\mathbb{Z}_2    & otherwise \\
\end{cases} $$ Thus for any simple square-cordial graph we have:
$$H_1(X_{(3)})= \Z^{p_0(G_{\Delta})^{coh}} \oplus 
\Z_2^{p_0(G_{\Delta})-p_0(G_{\Delta})^{coh}}.$$
By Proposition 4.6,
 $H_1(X_{\Delta}) = \Z^{p_0(G_{\Delta})^{coh}} \oplus
\Z_2^{p_0(G_{\Delta})-p_0(G_{\Delta})^{coh}} \oplus 
\textrm{tor}_3H_1(X_{\Delta})$, 
and $\textrm{tor}_3H_1(X_{\Delta})$ has no $\Z_9$ torsion.\\ 

Furthermore, we have $H_1(X_3)= H_1(\hat X_3)\oplus \Z^{v-1}$. 
Combining these arguments together and applying Proposition 4.6, we get \\
$\textrm{tor}_3H_1(X_{\Delta}) = 
\Z_3^{v-1 + dim H_1(\hat X_3,\Z_3)- p_0(G_{\Delta})^{coh}}$
as needed.
\end{proof}
\begin{corollary}\label{Corollary 4.8}
 If $G$ is a simple square-cordial graph then 
\begin{eqnarray*}
&& H^{1,2v-3}_{\A_3}(G) \\ &=& \Z_2^{p_0(G_{\Delta})-p_0(G_{\Delta})^{coh}} \oplus
\Z_3^{v-1 + dim H_1(\hat X_3,\Z_3)- p_0(G_{\Delta})^{coh}} \oplus 
\Z^{p_0(G_{\Delta})^{coh} + t_2 + 2t_3 - E - \frac{d_2}{2}}.
\end{eqnarray*}

\end{corollary}
Note that $H_1(X_{\Delta}(G))= H_1(X_{\Delta,4}(G))$ by Remark 4.4 -- 
we do not have to add 2-cells along 4-cycles in square-cordial graphs.
Now proof follows directly from Theorem 4.2(0) and Lemma 4.7 .

\begin{corollary}\label{Corollary 4.9}
For the complete graph with $n$ vertices $K_n$, $n\geq 4$ we have
$$H^{1,2n-3}_{\A_3}(K_n)= 
\mathbb{Z}_2\oplus \mathbb{Z}_3^{n-1}\oplus\Z^{\frac{n(n-1)(2n-7)}{6}}.$$
\end{corollary}
\begin{proof} Since $K_n$ is a square-cordial graph, $G_{\Delta}$ is connected.
 Consequently,  $H_1(X_{3}) = 0$ (as every cycle 
is a boundary cycle) and for $n>3$, $G_{\Delta}$ is not coherent.  
Therefore, by Lemma 4.7: 
$H_1(X_{\Delta}(K_n))= \Z_2 \oplus \Z_3^{n-1}$.
For the complete graph we have $t_0=d_2=d_{\geq3}=0$, so from 
Theorem 4.2(0) we get $H_{0,2n-3}^{\A_3}(K_n)= \Z_2 \oplus \Z_3^{n-1}$. 
Because $H_{0,2n-3}^{\A_3}(K_n)$ is a torsion group,
the chain map $C_{1,2v-3}^{\A_3} \to C_{0,2v-3}^{\A_3}$ is an 
epimorphism over $Q$. Thus 
\begin{eqnarray*}
&& rank\, H_{1,2n-3}^{\A_3}(K_n) = 
rank\, C_{1,2v-3}^{\A_3} - rank\, C_{0,2v-3}^{\A_3} \\
&=& \frac{n(n-1)}{2}(n-1) - 
\left( \frac{n(n-1)(n-2)}{6} + n(n-1)\right) = \frac{n(n-1)(2n-7)}{6}. 
\end{eqnarray*}
Therefore $H_{1,2v-3}^{\A_3}(K_n)= \Z^{\frac{n(n-1)(2n-7)}{6}}$ and \\
$H^{1,2n-3}_{\A_3}(K_n)= tor H_{0,2n-3}^{\A_3}(K_n) \oplus 
H_{1,2n-3}^{\A_3}(K_n) = 
\mathbb{Z}_2\oplus \mathbb{Z}_3^{n-1}\oplus\Z^{\frac{n(n-1)(2n-7)}{6}}.$
The proof of Corollary is completed\footnote{For the triangle $K_3$ we 
are getting $H_{0,3}^{\A_3}(K_3)=\Z \oplus \Z_3$, $H_{1,3}^{\A_3}(K_3)=0$ 
and $H^{1,3}_{\A_3}(K_3)=\Z_3$ which do not agree with the formula from 
Corollary 4.9 because the graph $G_{\Delta}$ is coherent.}.
\end{proof}
\begin{corollary}\label{Corollary 4.10}
If a simple graph $G$ contains a triangle then $H^{1,2v(G)-3}_{\A_3}(G)$
contains $\mathbb{Z}_3$.
\end{corollary}
\begin{proof}
If $G$ has an oriented triangle with edges $e_1,e_2,e_3$ we have
the following relation in $H^{\A_3}_{0,2v-3}(G)$:
$$ 3(e_1-e_2)=3(e_2-e_3)=3(e_3-e_1)=0.$$
We need to prove that $e_1-e_2$ is not $0$ in homology. Assume
that $e_1-e_2=0$. As $G$ embeds in the complete graph, $e_1$,
$e_2$ are also elements of $H^{\A_3}_{0,2v-3}(K_{v})$ and for
$K_{v}$ we have more relations then for $G$, thus the equation
$e_1-e_2=0$ is valid also in $H^{\A_3}_{0,2v-3}(K_{v})$. By
symmetry of the complete graph we get that all edges are equal in
homology group $H^{\A_3}_{0,2v-3}(K_{v})$. Thus the homology
cannot contain $\mathbb{Z}_3$. This contradicts 
Corollary 4.10 and thus $H^{\A_3}_{0,2v-3}(G)$ 
contains $\mathbb{Z}_3$ torsion.
\end{proof}

\begin{corollary}\label{Corollary 4.11}
Let $W_k$ denote the wheel, that is, the graph which is a cone over 
$(k-1)$-gon (Figure 8.3). Then \\
$$H^{1,2n-3}_{\A_3}(W_n)=\left\{
\begin{array}{ll}
    \mathbb{Z}_3^{n-2} \oplus \mathbb{Z}^n, & \hbox{if n odd;} \\
   \mathbb{Z}_3^{n-1} \oplus \mathbb{Z}_2 \oplus \Z^{n-1}, & \hbox{if
n even.} \\
\end{array}%
\right.$$
\end{corollary}

\begin{proof} Corollary 4.11 follows from Lemma 4.7 by observing that 
the wheel $W_n$ is a square cordial-graph 
with $H_1(X_3)=0$ and $(W_n)_{\Delta}$ is connected and 
coherent if and only if $n$ is odd.
\end{proof}
In the next section we describe initially unexpected examples
of torsion in cohomology, finding in particular that 
for any $k$ there is a graph $G$ such that $H^{1,2v-3}_{\A_3}(G)$ 
contains $Z_k$. Graphs we consider are not square-diagonal; nevertheless 
Lemma 4.7 is very useful in the analysis.

\section{The family of $G_n$ graphs}\label{Section 5}
The following graph, denoted by $G_4$, has $\Z_5$ in cohomology,
that is, \\ $H^{1,37}_{\A_3}(G_4)=
\Z_2 \oplus \Z_3^{19} \oplus \Z_5 \oplus \Z^{25}$.\ \\

\centerline{\psfig{figure=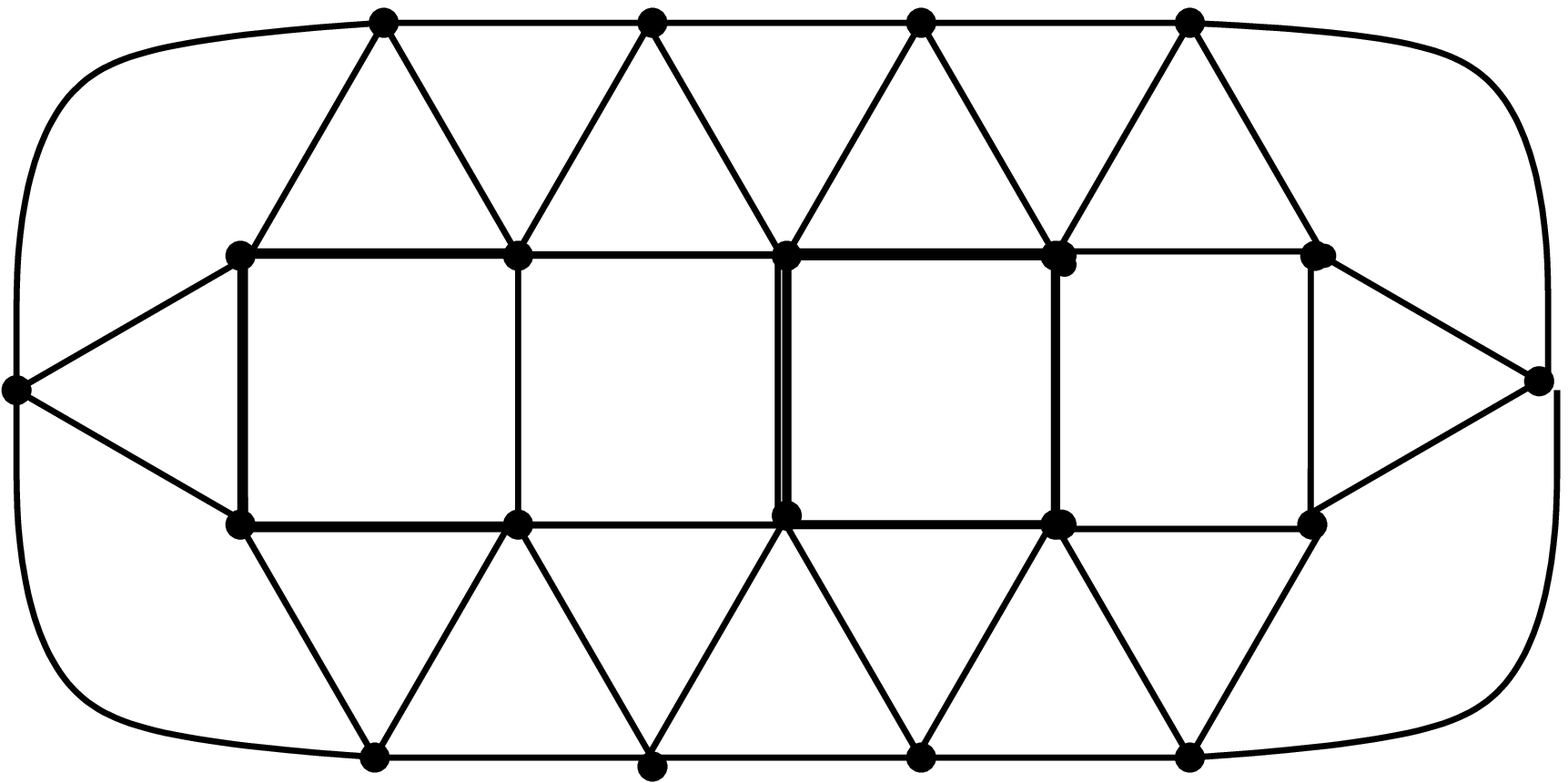,height=4.9cm}}

\centerline{Figure 5.1}
\ \\

Initially we were surprised to see torsion different from $\Z_3$ or $\Z_2$.
The reason was that 
A.~Shumakovitch has conjectured that for alternating links the
torsion in Khovanov homology can have only elements of order $2$.
We have conjectured analogously that for the algebra $\A_2$ the torsion
part of $H^{*,*}_{\A_2}(G)$ can  have only elements of order $2$.
This is still an open problem, and motivated by this we thought 
that for $\A_3$ graph homology torsion will be rather limited.
However, after computing  $H^{1,5}_{A_3}(K_4) = \Z_3^2 \oplus \Z_6
\oplus \Z^2$ we showed that torsion can have elements of different order. 
Furthermore, we
analyzed the family of plane graphs $G_1$, $G_2$,...,$G_k$ (see
Figures 5.1, 5.2) and we had a strong indication that for
any $n$ there is a graph with torsion $\Z_n$. In particular,
$H^{1,13}_{\A_3}(G_1)= \Z_4 \oplus \Z_3^7 \oplus \Z^{12}$,
$H^{1,21}_{\A_3}(G_2)= \Z_{18}\oplus \Z_3^{10}\oplus \Z^{15}$,
$H^{1,29}_{\A_3}(G_3)= \Z_{8}\oplus \Z_3^{15}\oplus \Z^{20}$,
$H^{1,37}_{\A_3}(G_4)= \Z_{10}\oplus \Z_3^{19}\oplus \Z^{25}$,
$H^{1,45}_{\A_3}(G_5)= \Z_{4}\oplus \Z_9 \oplus \Z_3^{22}\oplus
\Z^{30}$, $H^{1,53}_{\A_3}(G_6)= \Z_{14}\oplus \Z_3^{27}\oplus
\Z^{35}$, and $H^{1,61}_{\A_3}(G_7)= \Z_{16}\oplus
\Z_3^{31}\oplus \Z^{40}$.\\
Graphs $G_n$ are not square-diagonal so we cannot use Lemma 4.7 
directly, however combining it with 
  the Main Lemma we have
proved that $H^{1,8(k+1)-3}_{\A_3}(G_k)$ contains $\Z_{6k+6}$
or more generally:
\begin{corollary}\label{Corollary 5.1}
For $k>1$ we have:
$$H^{1,8(k+1)-3}_{\A_3}(G_k)= \Z_{6k+6} \oplus \Z_3^{4k+2}
\oplus \Z^{5k+5}$$
\end{corollary}
\begin{proof}
Let $G'_k$ be the graph obtained from $G_k$ by deleting all edges
which are not on 3-cycles (say $f_1,f_2,...,f_{k-1}$).
$G'_k$ is a square-cordial graph with the 
triangle graph $(G'_k)_{\Delta}$ connected and coherent. Therefore from
Lemma 4.7 it follows that
 $H_1(X_{\Delta}(G'_k))=Z^{E(G'_k)}/I_{\Delta} =Z \oplus \Z_3^{4k+3}$.
Furthermore,  $H_1(X_{\Delta,4}(G))=
\Z^{E(G)}/I_{\Delta,4} = \Z^{E(G'_k)}/I_{\Delta,4} \cup R$, where
$R$ is the relation in which the sum of the edges along rectangle
composed of squares is equal to 0. Combining this with previous
relations we get the relation $(2k+2)e + h_{\Delta}$, where $e$ is 
a generator of $\Z$ and $h_{\Delta}$ belongs to $\Z_3^{4k+3}$. 
The 3-torsion element $h_{\Delta}$ cannot be equal to zero because 
$$(\Z^{E(G)}/I_{\Delta,4}) \otimes \Z_3 = 
H_1(X_{3,4}(G),\Z_3) =  \Z_3^{4k+3} \oplus H_1(\hat X_{3,4}(G),\Z_3) = 
\Z_3^{4k+3}.$$ Therefore $H_1(X_{\Delta,4}(G))= 
\Z_{6k+6} \oplus \Z_3^{4k+2}$. To find the free
part of $H^{1,8(k+1)-3}_{\A_3}(G_k)$ we use the standard tricks
using Euler characteristic and duality (homology -- cohomology).
\end{proof}
\ \\
\ \\
\begin{figure}[h]
\begin{center}
\scalebox{.3}{\includegraphics{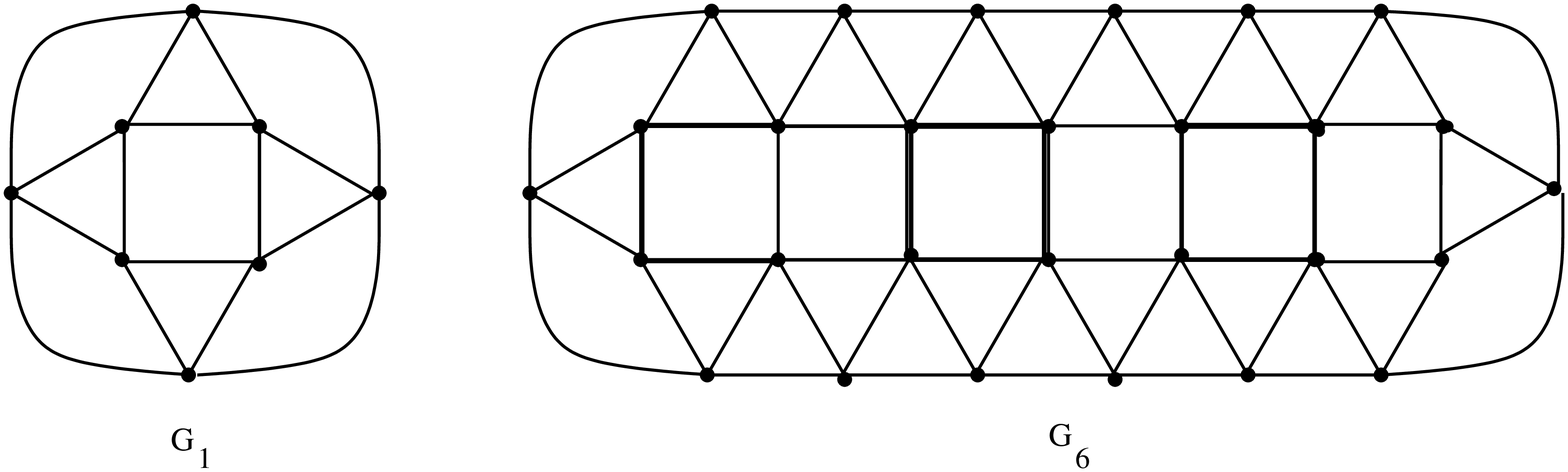}}
\end{center}
\end{figure}\ \\
\centerline{Figure 5.2;\ Family of graphs $\{G_i\}$, where $i$ is the number
of squares}
\ \\ \

We have constructed in this section, for any $n$,
  graphs which have $\Z_n$ in torsion of cohomology. In 
Corollary 4.10 we have proven that a simple graph with a 3-cycle 
has $\Z_3$ torsion in cohomology. However Theorem 4.1 allows 
us to construct any torsion, even for a graph with no 3-cycle.

\begin{corollary}\label{Corollary 5.2}
For any finitely generated abelian group $T$ there is a simple 
graph, $G$, without any 3-cycle 
such that $tor H^{1,2v-3}_{\A_3}(G) = T$.
\end{corollary}
\begin{proof} 
For a simple graph, $G$ without a 3-cycle, the torsion of 
$H^{1,2v-3}_{\A_3}(G)$ is equal to the torsion of $H_1(X_4(G)$ where, 
$X_4(G)$ is obtained from $G$ by attaching a 2-cell along every 4-cycle.
It is not difficult to construct $G$ such that 
$tor X_4(G) = T$. For example, to obtain $T=\Z_2$ we divide a projective 
plane into sufficiently many squares (25 suffices). In Figure 5.3 
we present a graph with no 3-cycles and $T=\Z_3$ (square divided into 
$6$ by $6$ small squares and boundary quotiented by $\Z_3$ action).
\end{proof}

\centerline{\psfig{figure=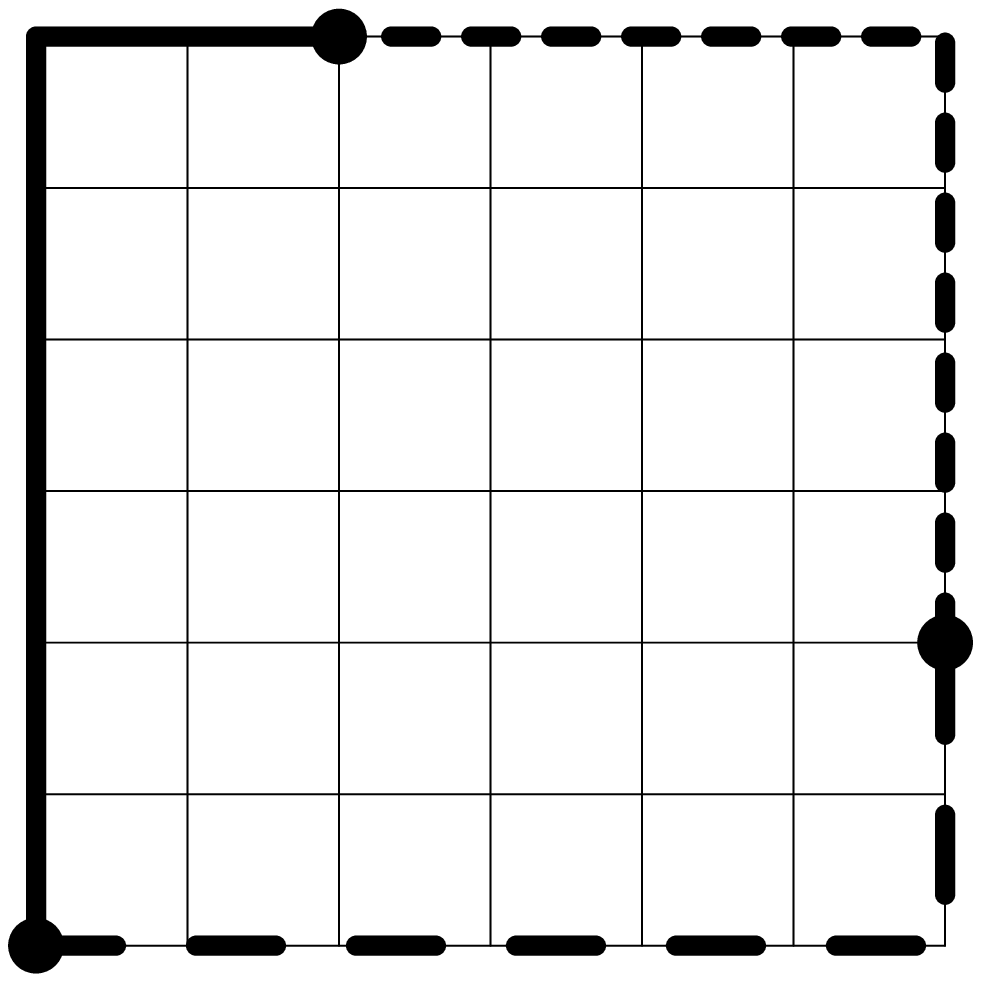,height=3.8cm}}

\centerline{Figure 5.3: A graph with $tor H^{1,63}_{\A_3}(G) = \Z_3$}

\section{Examples of graphs with the same dichromatic (and Tutte) polynomial
but different $\A_3$ first cohomology}
In this section we construct graphs which have the same chromatic (even
dichromatic and Tutte)
polynomial but different $H^{1,2v-3}_{\A_3}(G)$. In relation to these
examples we also describe $A_3$ graph cohomology for the one vertex 
product of graphs. We work also with the 
2-vertex product of graphs, showing in 
particular that the Whitney flip (reglueing of vertices) preserves the  
cohomology of the product\footnote{Whitney flip of graphs is closely 
related to mutation of links.}. In effect, 
$H^{1,2v-3}_{\A_3}(G)$ is a 2-isomorphism (i.e. matroid) graph invariant.

One of the corollaries of our main theorem (Theorem 4.1) is that we
can give simple formulas for cohomology of vertex 
products of graphs, $G*H$ and and edge
products of graphs,  $G|H$, (if working with $\Z_3$ coefficients),
 where $G|H$ is obtained from
$G$ and $H$ by identifying an edge in $G$ with an edge in $H$.
An edge product depends on
the choice of identified edges (see Fig. 6.1 for simple examples) but
the chromatic polynomial ($P(G)\in \Z[\lambda]$) is always the same,
$P(G|H)=P(G)P(H)/(\lambda  (\lambda -1))$ (see e.g. \cite{Big}).
However we can often differentiate between different products $G|H$
using $H^{1,2v-3}_{\A_3}(G|H)$, for example for the first pair in
Figure 6.1 we get, $\Z_3^2 \oplus \Z$ and $\Z_3^2$, respectively and
for the second pair in
Figure 6.1 we get, $\Z_3^3 \oplus \Z^4$ and $\Z_3^3 \oplus \Z^3$, 
respectively.\\
\ \\

\ \\
\centerline{\psfig{figure=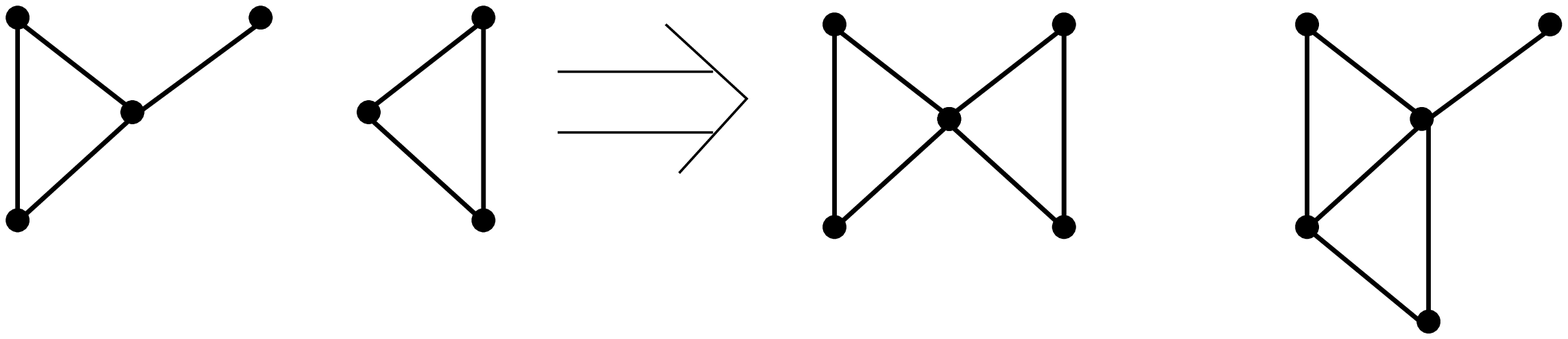,height=2.8cm}}\ \\
\centerline{\psfig{figure=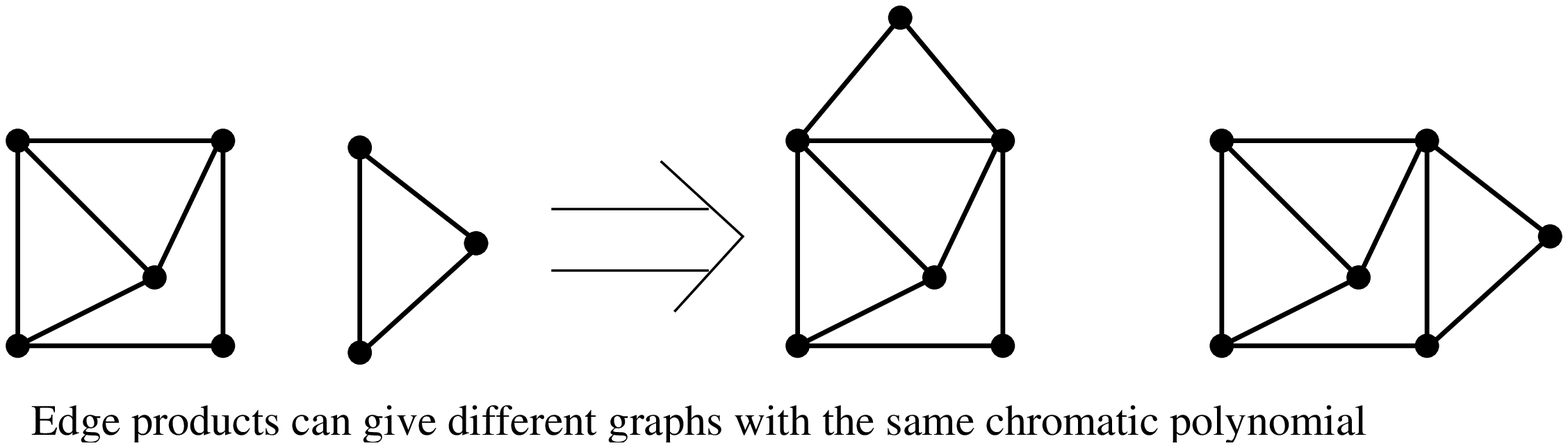,height=3.8cm}}
\centerline{ Figure 6.1: }
\ \\

In the 1930's Marion Cameron Gray found the following example of
graphs which are not
2-isomorphic but have the same dichromatic (so also Tutte)
polynomials, Figure 6.2, \cite{Tut,Big}. These examples have different first
$\A_3$ graph cohomology. In fact they differ from the second
example in Figure 6.1 only by multiple edges, thus
$H^{1,9}_{\A_3}(\textrm{Gray}_1)= \Z^3_3 \oplus Z^4$
and $H^{1,9}_{\A_3}(\textrm{Gray}_2)= \Z^3_3 \oplus Z^3$.
\ \\ \ \\
\centerline{\psfig{figure=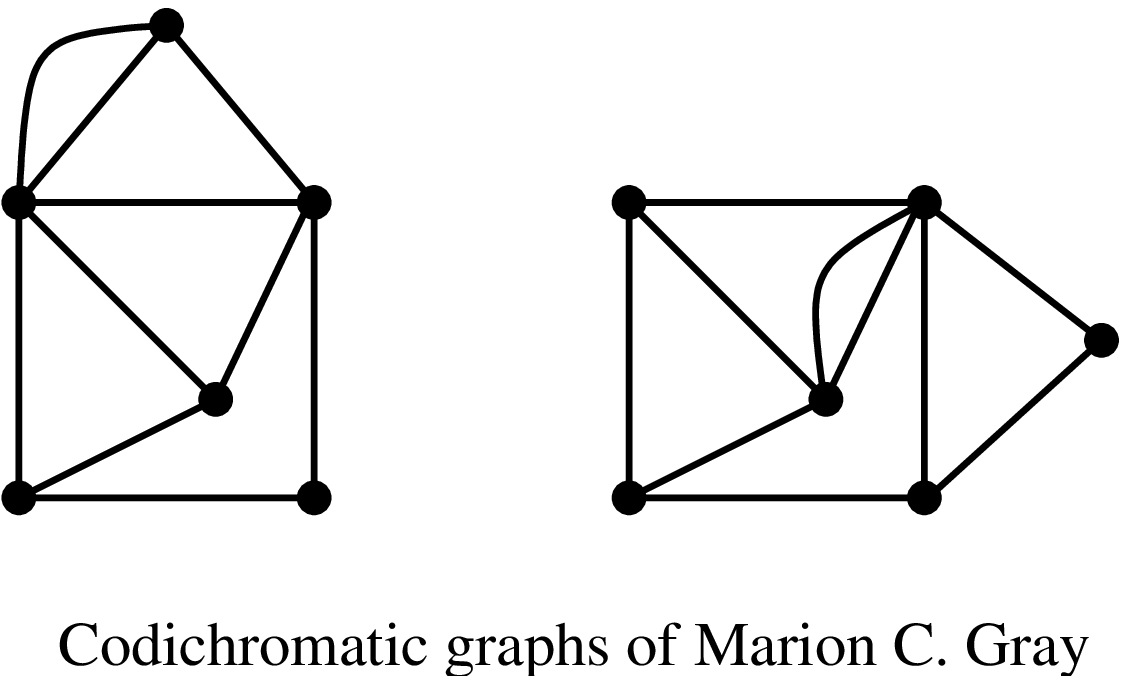,height=4.9cm}}
\centerline{ Figure 6.2: }
\ \\ \

The Gray graphs are not 2-isomorphic (see Definition 6.1) however as 
an example of codichromatic graphs with different $\A_3$ graph cohomology 
they are a little disappointing as they have multiple edges. 
In the sequel paper we will analyze examples obtained 
by rotation \cite{Tut,APR}.

 We devote the 
rest of this section to the proof that 2-isomorphic graphs have isomorphic 
cohomology $H^{1,2v-3}_{\A_3}$. First recall definition of 2-isomorphism 
in the form convenient for our considerations, i.e. based on one-vertex, 
$G*H$, and two-vertex, $G^*_*H$, products.

\begin{definition}\label{Definition 6.1}
\begin{enumerate}
\item[(1)] A one-vertex product $G*H =G*(v=w)H$ of graphs $G$ and $H$ 
with base points $v$ and $w$ respectively is obtained by gluing 
$G$ with $H$ by identifying $v$ and $w$.
\item[(2)]
We define a two-vertex product 
$G^*_*H =G_{*(v_2=w_2)}^{*(v_1=w_1)}H$ as follows.
Given a graph $G$ with two chosen vertices $v_1, v_2$ and a graph $H$ 
with chosen vertices $w_1$ and $w_2$ we obtain a two-vertex product of 
$G$ and $H$ by identifying $v_1$ with $w_1$ and $v_2$ with $w_2$. 
If we switch the roles of $w_1$ and $w_2$ in a two-vertex product we 
obtain the new graph $G_{*(v_2=w_1)}^{*(v_1=w_2)}H$ (in standard terminology 
we say that these two graphs differ by Whitney flip).
\item[(3)]
2-isomorphism of graphs is the smallest equivalence relation on isomorphism 
classes of graphs which satisfies:\\
(i) One-vertex products of $G$ and $H$ (for any attachments) are 
2-isomorphic.\\
(ii) Two graphs which differ by Whitney flip are 2-isomorphic.
\end{enumerate}
\end{definition}

\begin{theorem}\label{Theorem 6.2}\ \\  
(1) Let $v$ and $w$ be two vertices on the graph $Y$ 
(not necessary connected), and let $Y'=Y/(v=w)$ 
be the graph obtained from $Y$ by identifying $v$ with $w$. If the distance 
$d_Y(v,w)> 4$ (allowing $d_Y(v,w)= \infty$), 
then $H^{1,2v(Y')-3}_{\A_3}(Y') = H^{1,2v(Y)-3}_{\A_3}(Y)$.
In particular, $$H^{1,2v(G*H)-3}_{\A_3}(G*H) = 
H^{1,2v(G)-3}_{\A_3}(G) \oplus H^{1,2v(H)-3}_{\A_3}(H).$$
(2) 
$H^{1,2v(G_{*}^{*}H)-3}_{\A_3}(G_{*(v_2=w_2)}^{*(v_1=w_1)}H)
=H^{1,2v(G_{*}^{*}H)-3}_{\A_3}(G_{*(v_2=w_1)}^{*(v_1=w_2)}H)$.\\
(3) 
 If two graphs, $G_1$ and $G_2$ are 2-isomorphic
then $H^{1,2v-3}_{\A_3}(G_1)= H^{1,2v-3}_{\A_3}(G_2)$.
\end{theorem}
\begin{proof}
Without loss of generality we can assume that we deal only with simple 
graphs (in case when our operation produce a multiple edge we will 
replace  it by a singular edge without affecting cohomology).\\
(1) It follows directly from Theorem 4.1: first the cell complex 
$X_{\Delta,3}(Y')$ is equal to $X_{\Delta,3}(Y)$. Furthermore, 
$t_2 - \frac{d_2}{2}$ and $sq$ are the same for $Y'$ and $Y$.\\
(2) We will consider three case depending on position of $v_1$, $v_2$, $w_1$ 
and $w_2$ in respective graphs.
(i) We assume that either $v_1$ is connected by an edge to $v_2$ in $G$ 
or $w_1$ is connected by an edge to $w_2$ in $H$. \\
(ii) Assume that (i) does not hold and the distance 
$$d_G(v_1,v_2) + d_H(w_1,w_2) > 4.$$
(iii) Assume that neither (i) nor (ii) hold, that is, 
$$d(v_1,v_2)=2=d(w_1,w_2).$$
First notice that 
\begin{eqnarray*}
&&
(t_2-\frac{d_2}{2} -sq)(G^*_*H)\\
& =& 
(t_2-\frac{d_2}{2} -sq)(G) + 
(t_2-\frac{d_2}{2} -sq)(H) + \left\{
\begin{array}{ll}
0& \textrm{in the cases  (i), (ii)}\\
1& \textrm{in the case (iii)}
\end{array}
\right.
\end{eqnarray*}
To complete our proof of the theorem we should show, according to 
Theorem 4.1, that $H^2(X_{\Delta,4},\Z)(G^*_*H)$ is invariant under 
Whitney flip. In the sequel paper we will discuss precise formulas 
relating $H^2(X_{\Delta,4}(G^*_*H),\Z)$ with $H^2(X_{\Delta,4},\Z)(G)$ 
and $H^2(X_{\Delta,4},\Z)(G)$, compare Proposition 6.3 ; here we only 
notice that Whitney flip of $G^*_*H$ changes the chain complex 
of cell-complex $X_{\Delta,4}(G^*_*H))$ only by isomorphism.
We discuss concisely the cases (i), (ii) and (iii) separately.\\
(i) 
In this case we can think of  $G_{*(v_2=w_2)}^{*(v_1=w_1)}H$ as being
of the form
$G|H$ were $v_1, v_2$ are connected by an edge, 
$\overrightarrow{e_{G}}$, in  $G$ and
$w_1, w_2$ are connected by an edge, 
$\overrightarrow{e_{H}}$, in  $H$. Then 
we check that constructing $G|H$ by glueing
$\overrightarrow{e_{G}}$ to $\overrightarrow{e_{H}}$
or $\overrightarrow{e_{G}}$ to $\overleftarrow{e_{H}}$ 
yields the same cohomology. That is, chain complexes of 
$X_{\Delta,4}(G|(\overrightarrow{e_{G}}=\overrightarrow{e_{H}}H)$ 
and of $X_{\Delta,4}(G|(\overrightarrow{e_{G}}=\overleftarrow{e_{H}})$ 
are isomorphic (the isomorphism is sending $\overrightarrow{e}$ to 
$\overrightarrow{e}$ if $\overrightarrow{e}$ is an edge of $G$ and 
it is sending $\overrightarrow{e}$ to $\overleftarrow{e}$ if 
$\overrightarrow{e}$ is an edge of $H$. \\
(ii) The distance $d(v_2,w_2)$ in $G*(v_1=w_1)H$ is greater than $4$.
Therefore we can use (1) to conclude that
$$H^{1,2v(G^*_*H)-3}_{\A_3}(G^*_*H) = H^{1,2v(G)-3}_{\A_3}(G) \oplus
H^{1,2v(H)-3}_{\A_3}(H).$$
(iii) In this case $d(v_1,v_2)=2=d(w_1,w_2)$, therefore
$v_1$ is connected with $v_2$ in $G$ by a ``joint" $(e^{(1)}_G,e^{(2)}_G)$ 
in $G$ and by a ``joint" $(e^{(1)}_H,e^{(2)}_H)$ in $H$. 
As in the case (i) we can see that
$X_{\Delta,4}(G_{*(v_2=w_2)}^{*(v_1=w_1)}H)$ and 
$X_{\Delta,4}(G_{*(v_2=w_1)}^{*(v_1=w_2)}H)$ yield isomorphic 
chain complexes of cell-complex homology. As in (i) we 
send $\overrightarrow{e}$ to
$\overrightarrow{e}$ if $\overrightarrow{e}$ is an edge of $G$ and
we send $\overrightarrow{e}$ to $\overleftarrow{e}$ if
$\overrightarrow{e}$ is an edge of $H$.\\
(3) It follows from (1) and (2) by definition of 2-isomorphism
\end{proof}

As mentioned before, we can use Theorems 4.1 and 4.2 to find the formulas 
for cohomology of an edge and 2-vertex products of graphs. We discuss these 
in a sequel paper, listing here three easy but useful special cases.
\begin{proposition}\label{Proposition 6.3}\
\begin{enumerate}
\item[(i)]
Let $G$ and $H$ be simple graphs and\\
 $G|H = G|(\vec e_G=\vec e_H)H$, where $\vec e_G$, $\vec e_H$ are identified edges. Then
$$H^{1,2v(G|H)-3}_{\A_3}(G|H,\Z_3) = H^{1,2v(G)-3}_{\A_3}(G,\Z_3) \oplus
H^{1,2v(H)-3}_{\A_3}(H,\Z_3) \oplus \Z_3^{t(e_G)t(e_H)},$$
where
$t(e_G)$ is the number of 3-cycles in $G$ containing $e_G$.

\item[(ii)] $H^{1,2v(G)-1}_{\A_3}(G|P_3) =
H^{1,2v(G)-3}_{\A_3}(G) \oplus \Z_3 \oplus \Z^{t(e_G)}$.
\item[(iii)] $H^{1,2v(G)}_{\A_3}(G|P_4) =
H^{1,2v(G)-3}_{\A_3}(G) \oplus \Z$.
\end{enumerate}
\end{proposition}
\ \\
Sketch of a proof of (i).\\
First we find $H_1(X_{\Delta,4}(G|H,\Z_3))$. As noted in Remark 4.4,
$H_1(X_{\Delta,4}(G|H))=H_1(X_{\Delta,4'}(G|H))$ and\footnote{For 
the convenience of a reader we recall useful definitions 
of various cell-complexes built from $G$:\\
1(i) $\hat X_3(G)$ is the cell complex obtained from $G$ by 
attaching 2-cells along 3-cycles in $G$,\\
1(ii) $\hat X_{3,4'}(G)$ is the cell complex obtained from $\hat X_3(G)$ by 
attaching 2-cells along 4-cycles with no diagonals.\\
1(iii) $\hat X_{3,4'}(G)$ is the cell complex obtained from $\hat X_3(G)$ by
attaching 2-cells along 4-cycles.\\
2(i) $\hat X_{(3)}(G)$ is the cell complex obtained from $G$ by
 identifying edges of every triangle
in a coherent way (i.e. $\vec e_i = \vec e_j =\vec e_k$ for a triangle 
oriented as in Figure 4.8).\\
2(ii) $\hat X_{(3),4'}(G)$ is the cell complex obtained from $\hat X_{(3)}(G)$ 
by attaching 2-cells along 4-cycles with no diagonals.\\
2(iii) $\hat X_{(3),4}(G)$ is the cell complex obtained from $\hat X_{(3)}(G)$
by attaching 2-cells along 4-cycles.\\
3. If $\hat X(G)$ denotes one of cell complexes defined in 1 or 2 the 
$ X(G)$ is obtained from $\hat X(G)$ by identifying all vertices of $G$.\\
4(i) $ X_{\Delta}(G)$ is a cell complex obtained from $G$ by 
identifying vertices of $G$ and then adding 2-cells
along expressions $2\vec e_3 - \vec e_2 - \vec e_1$
for any 3-cycle in $G$ -- two 2-cells added per 
every 3-cycle (see Fig. 4.7); compare Footnote 7.\\
4(ii) $ X_{\Delta,4'}(G)$ is the cell complex obtained from $ X_{\Delta}(G)$ 
by attaching 2-cells along 4-cycles with no diagonals.\\
4(iii) $ X_{\Delta,4}(G)$ is the cell complex obtained from $ X_{\Delta}(G)$
by attaching 2-cells along 4-cycles.} 
in $X_{\Delta,4'}$
every 2-cell is glued either to $G$ or to $H$. We know, Remark 4,4, that
$H_1(X_{\Delta,4}(G|H))=H_1(X_{\Delta,4'}(G|H))$ and \\
$H_2(X_{\Delta,4}(G|H)) = H_2(X_{\Delta,4'}(G|H)) + \Z^{sq'(G|H)}$.
Furthermore, if we work over $\Z_3$ we have 
$H_1(X_{\Delta,4'}(G|H),\Z_3)= H_1(X_{3,4'}(G|H),\Z_3)$ and\\ 
$H_2(X_{\Delta,4'}(G|H),\Z_3)= H_2(X_{3,4'}(G|H),\Z_3) \oplus \Z^{t_3}$. 
Since $e_G=e_H$ is contractible in $\hat X_{3,4'}(G|H)$ we have 
$H_2(\hat X_{3,4'}(G|H),\Z_3) =H_2(\hat X_{3,4'}(G),\Z_3)
 \oplus H_2(\hat X_{3,4'}(H),\Z_3)$.
These, taken together, give:\\
$H_2(X_{\Delta,4}(G|H),\Z_3)= H_2(X_{\Delta,4}(G),\Z_3) \oplus 
H_2(X_{\Delta,4}(H),\Z_3) + \Z_3^{t(e_G)t(e_H)}$\\
 as $sq'(G|H) = sq'(G) + sq'(H) + t(e_G)t(e_H)$. 

\begin{example} For graphs of Gray we have: 
$H^{1,9}_{\A_3}(Gray)= \Z^3_3 \oplus Z^4$. Each of these graphs is 
built as an edge product of the broken wheel $W_4^{in}$, Figure 8.3, triangle
and doubling two of the edges.
In homology we start from the group $\Z^2_3 \oplus Z^4$. 
Adding a triangle along an edge
 transforms it to $\Z^3_3 \oplus Z^4$, while 
doubling edges preserves it.  
\end{example}

\section{Dichromatic graph cohomology}\label{Section 7}
M.~Stosic observed that chromatic graph cohomology can be slightly
deformed to give a categorification of the 
dichromatic polynomial \cite{Sto}.
\begin{enumerate}
\item[(i)]
We define the dichromatic
 cohomology $\tilde H^{i,j}_{\A_m},\M(G)$ as follows:\\
 We consider the algebra of truncated polynomials
$\A_m= \Z[x]/(x^m)$
and $M= \A_m$ or $M= (x^{m-1})$.
We modify
$d_e$ (from the definition of the chromatic chain complex) to
get the chain complex $C^{*,*}_{\A_m,\M}(G)$. The only modification is
in the case $e$ has endpoints on the same component of
$[G:s]$. We put $\tilde d_e = x^{m-1}Id$. The grading which Stosic
associates to $x^i$ is equal to $m-1-i$ and the grading of an element
$(x^{i_0},x^{i_1},...,x^{k(s)-1})$ in $C^i_{\A_m}$ is given by
$i + k(s)(m-1) - \sum i_j$. Notice that for $i < \ell(G)$ we have
$(i + k(s) = v$.
In this way we have graded cochain complex and graded cohomology which can be
used to recover dichromatic polynomial (so it is  named {\it dichromatic
graph cohomology}).

\item[(ii)] It follows from the definition that for $i < \ell(G)$
we have  $d^i = \tilde d^i$. Therefore, for
$i < \ell(G)- 1$, $\tilde H^i_{\A_m}(G) = H^i_{\A_m}(G)$.
More precisely, $\tilde H^{i, v(m-1) - j}_{\A_m}(G) = H^{i,j}_{\A_m}(G)$.
Furthermore for $i = \ell(G)- 1$, we have: \\
$tor(\tilde H^{i, v(m-1) - j}_{\A_m}(G))= tor(H^{i,j}_{\A_m}(G))$.
\end{enumerate}

Most of our computations done of $H^{*,2v-3}_{\A_3}(G)$ is also valid 
for dichromatic graph cohomology (one should remember about grading 
transformation).

\section{Computational results, width of cohomology}\label{Section 8}
In previous sections, we completely computed 
$H_{\A_2}^{1,v-1}(G)$ and $H_{\A_3}^{1,2v-3}(G)$. Here we present a few
computational results and conjectures derived from them for $\A_m$ graph
cohomology, for $m\geq 3$.
In the first subsection we still work over $\A_3$ algebra searching for
gradings with non-trivial torsion. 
\par
\subsection{Width of $tor(H^{1,*}_{\A_3}(G))$}\ \\ \

Let $j_{max}$ be the maximal index $j$ such that
$tor\,H_{\A_m} ^{1,j}(G) \neq 0$ and
$j_{min}$ be the minimal index $j$ such that $tor\,H_{\A_m}^{1,j}(G) \neq 0$.
We define the width, $w^1_{\A_m}(G)$ of the torsion $tor\,H_{\A_m}^{1,*}(G)$
as follows:
\begin{displaymath}
w^1_{\A_m}(G)=
\left\{
\begin{array}{ll}
-1 & \textrm{if $H_{\A_m} ^{1,*}(G)$ are all free}\\
j_{max} - j_{min} &  \textrm{otherwise}
\end{array}
\right.
\end{displaymath}
As we mentioned before for an odd $n$, 
$torH^{1,*}_{\A_m}(P_n)= H^{1,\frac{n-1}{2}m}_{\A_m}(P_n)=\Z_m$
thus the width $w^1_{\A_m}(P_n)=0$. We also know that for any graph and algebra
$\A_3$, $j_{max}\leq 2v-3$ (Proposition 2.8). In fact, from 
results of Section 4 it follows that
for a simple graph $G$, $j_{max} = 2v-3$ iff $G$ contains a triangle
or $H_1(X_4(G))$ has a torsion. In this section we analyze $tor
H^{1,*}_{\A_3}(G)$ (in particular $j_{min}$) for several families of graphs
and conjecture a formula for width.

Let $P_{t,k}= P_3|P_{k}$ be a graph with $k+1$ vertices
 obtained by gluing a triangle 
and a $k$-gon along an edge, Figure 8.1. \ \\

\centerline{\psfig{figure=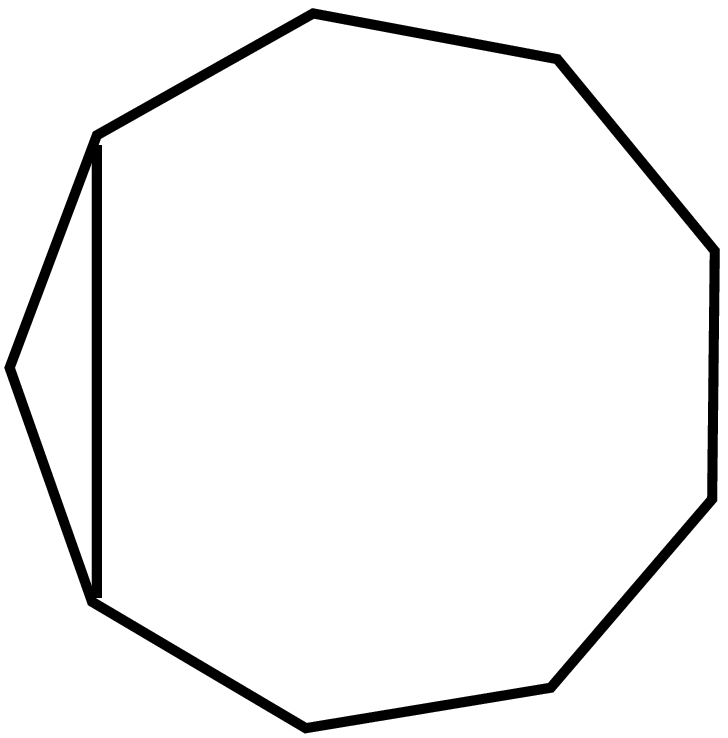,height=3.7cm}}
\centerline{ Figure 8.1:  $P_{t,8}= P_3|P_8$ } \ \\
\ \\

We know that 
$H^{1,2k-1}_{\A_3}(P_{t,k})=\Z_3$ for $k>4$ 
($H^{1,5}_{\A_3}(P_{t,3})=\Z_3^2\oplus \Z$ and 
$H^{1,7}_{\A_3}(P_{t,4})=\Z_3\oplus \Z$). Computational results are summarized
in the following Table and conjecture:

\begin{conjecture}\label{Conjecture 8.1}    
For graphs $P_{t,k}= P_3|P_{k}$ the following is true: $w^1_{\A_3}(P_{t,k})=k-3$                                 
\end{conjecture}

\begin{tabular}{|l|c| c| c| c| c| c|}\hline
 k & \textbf{$P_{t,k}$} &$torH_{\A_3} ^{1,2v-7}$& $torH_{\A_3} ^{1,2v-6}$& 
$torH_{\A_3} ^{1,2v-5}$& $torH_{\A_3} ^{1,2v-4}$& $H_{\A_3} ^{1,2v-3}$
 \\
&&&&&& \\ \hline 3 
&\psfig{figure=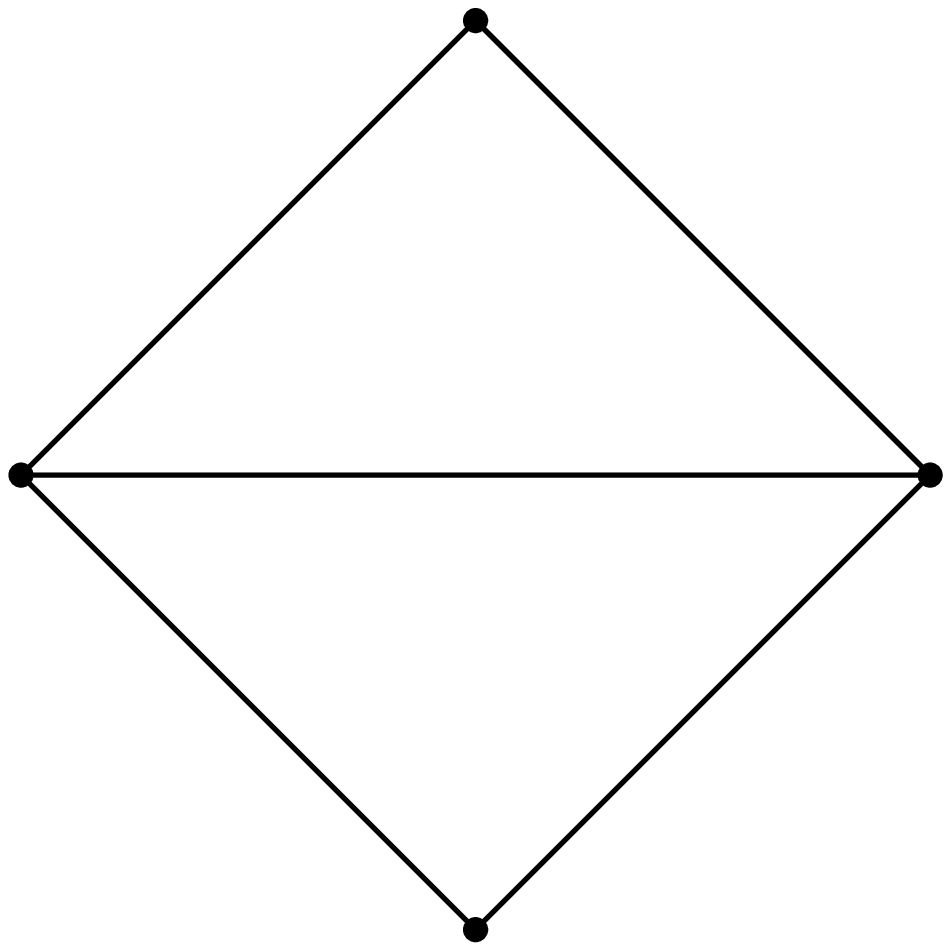,height=1.5cm}&0&0&0&0&$\mathbb{Z}
_3^2 \oplus \mathbb{Z}$  \\ \hline 4&\psfig{figure=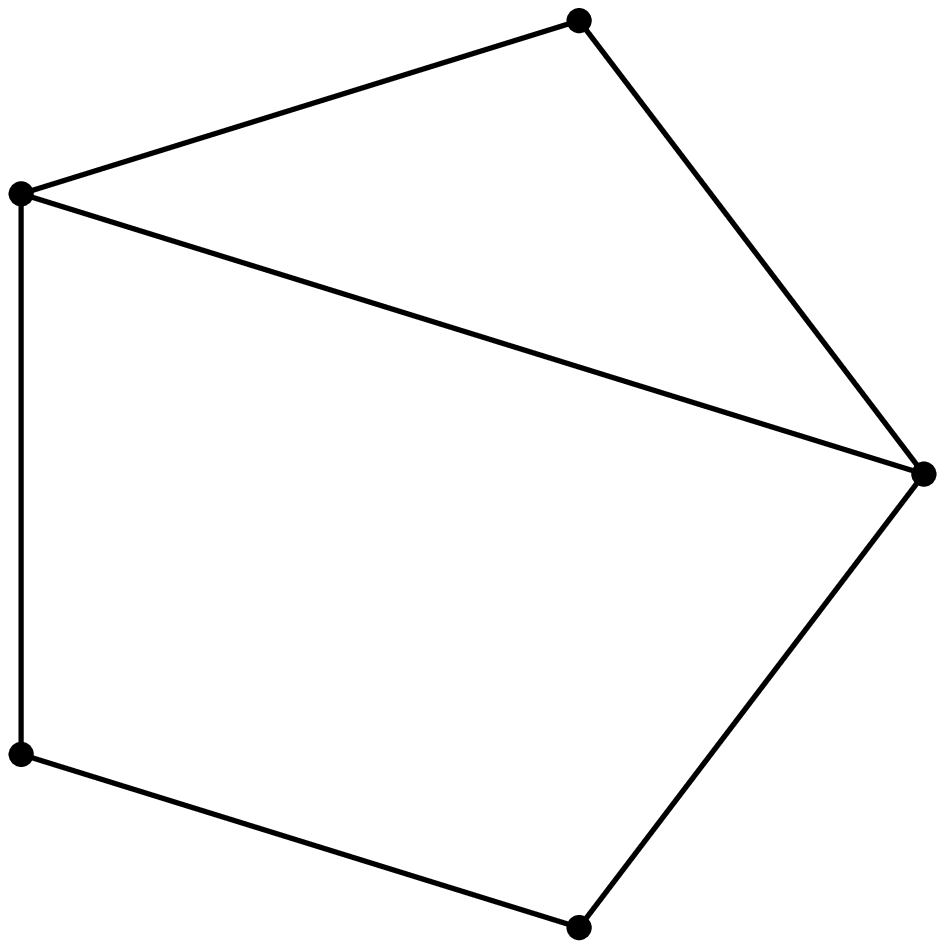,height=1.5cm} &
0&0&0&$\mathbb{Z}_3^2$&$\mathbb{Z}_3 \oplus \mathbb{Z}$ \\
\hline
 5&\psfig{figure=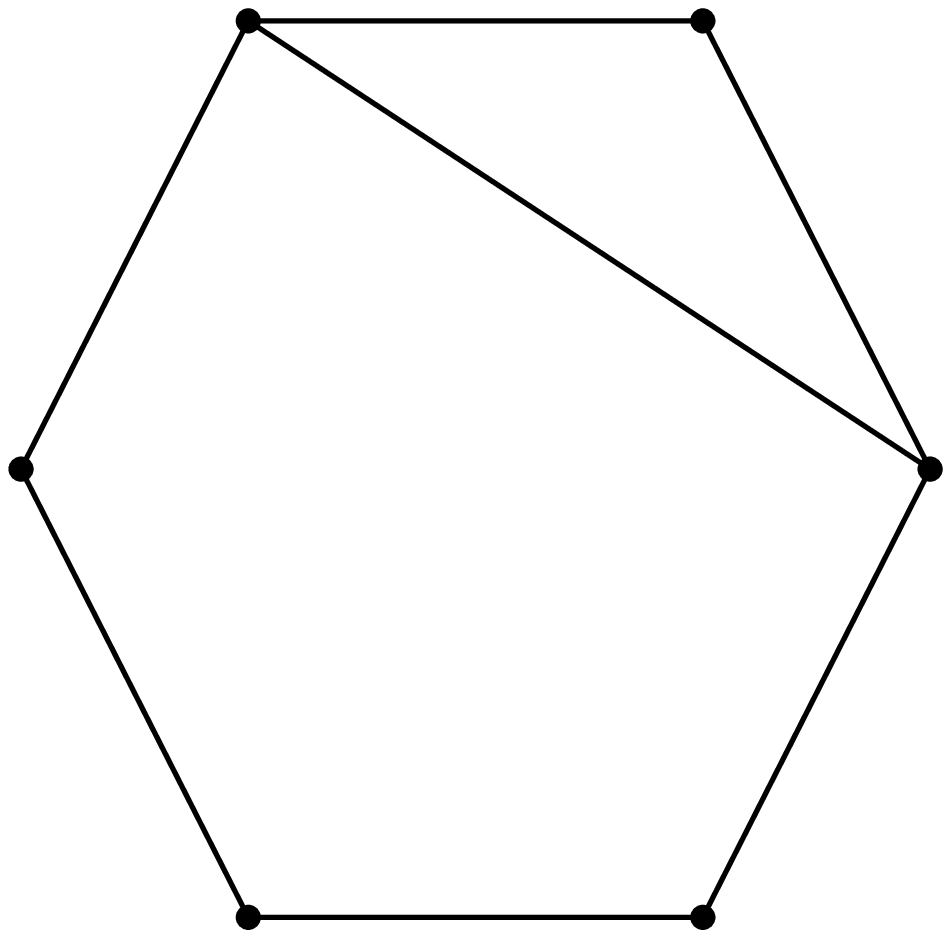,height=1.5cm} 
&0&0&$\mathbb{Z}_3^2$&$\mathbb{Z}_3
^{4}$&$\mathbb{Z}_3$ \\ \hline
6&\psfig{figure=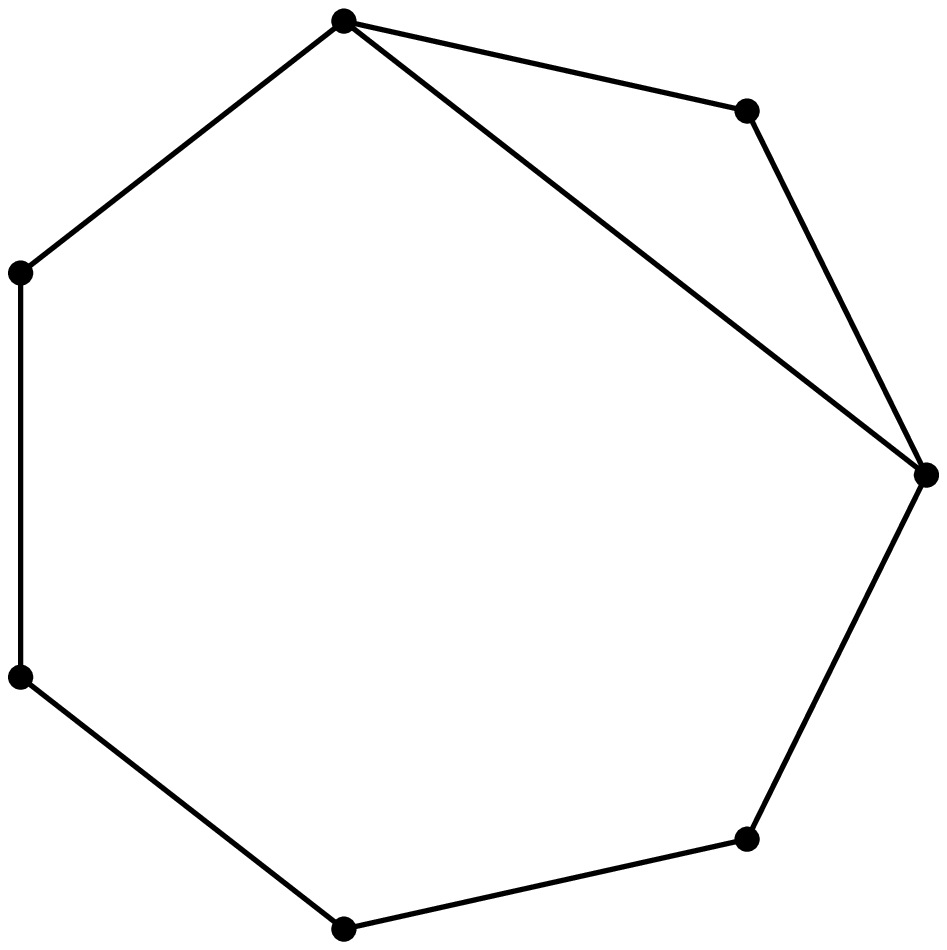,height=1.5cm}&0 
&$\mathbb{Z}_3$&$\mathbb{Z}_3^5$&$ \mathbb{Z}_3^4$&
$\mathbb{Z}_3$ \\ \hline

7&\psfig{figure=Triangle6.eps,height=1.5cm}&
$\mathbb{Z}_3$&$\mathbb{Z}_3^6$&$\mathbb{Z}_3^{10}$&$\mathbb{Z}_3^5$&
$\mathbb{Z}_3$ \\ \hline
8&\psfig{figure=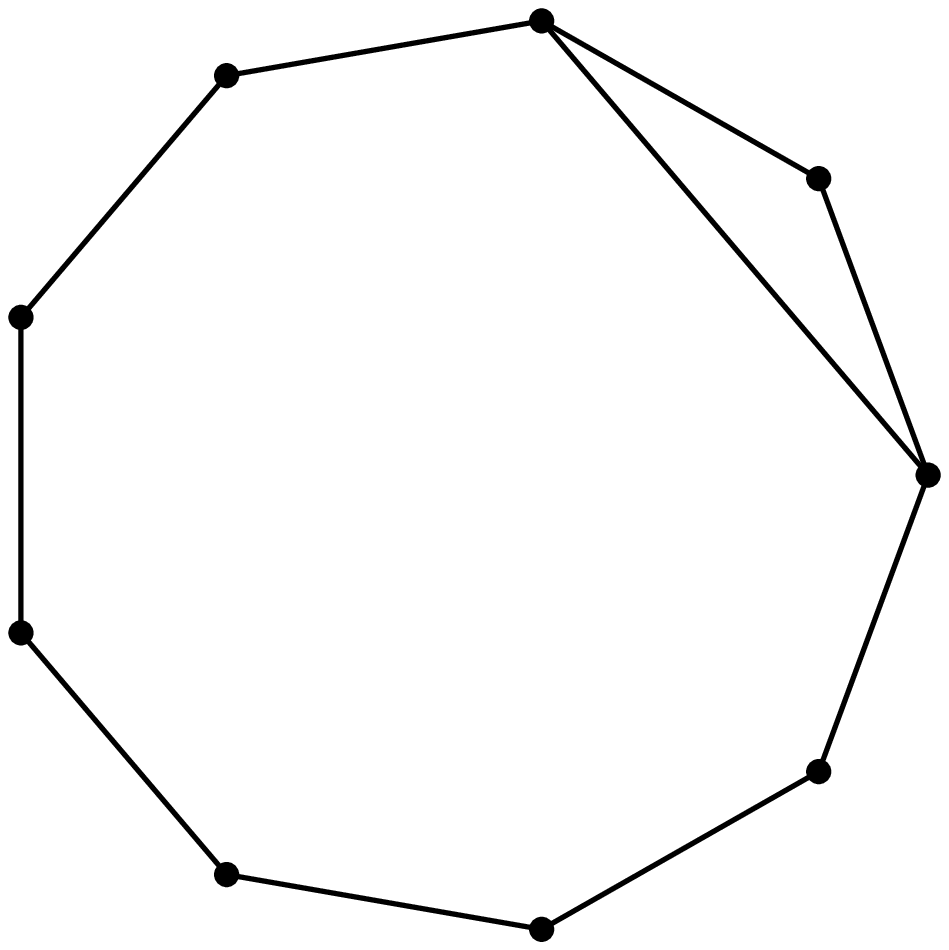,height=1.5cm}&$ 
\mathbb{Z}_3^6$&$\mathbb{Z}_3^{15}$
&$\mathbb{Z}_3^{14}$&$\mathbb{Z}_3^6$&
$\mathbb{Z}_3$ \\ \hline
\end{tabular}
\ \\
\ \\

 Another family of graphs we consider, $G_{t,s^k}$, is obtained by glueing a
 triangle and $k$ squares ($v=3+2k$ vertices)in a sequence along one edge;
 Figure 8.2.

\centerline{\psfig{figure=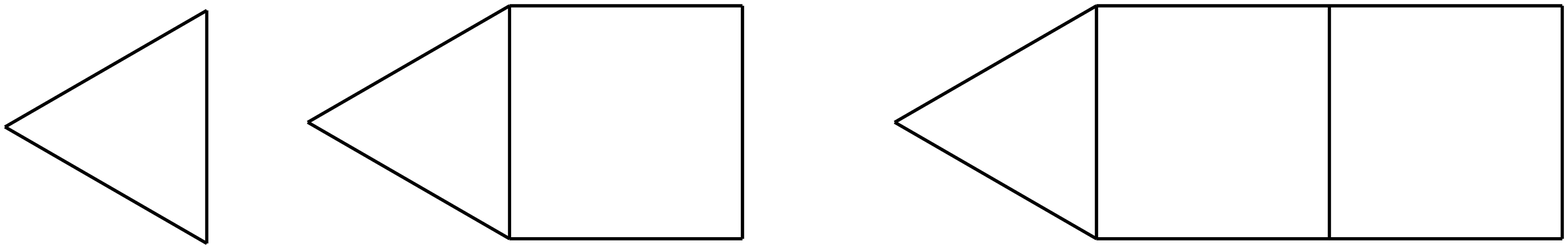,height=1.7cm}}
\centerline{ Figure 8.2: Family of $G_{t,s^k}$ graphs:\
$G_{t,s^0}=P_3$, $G_{t,s}$ and $G_{t,s^2}$} \ \\

From the theory developed earlier we get:
$H^{1,4k+3}_{\A_3}(G_{t,s^k})=\Z_3 \oplus \Z^k$ for $k\geq 0$.
Computational results are presented in the following conjectures and Table.

\begin{conjecture}\label{Conjecture 8.2}
$tor H^{1,4k+2}_{\A_3}(G_{t,s^k})= \mathbb{Z}_3^{2k}.$
\end{conjecture}
\begin{conjecture}\label{Conjecture 8.3}
$tor H^{1,3k+3}_{\A_3}(G_{t,s^k})= \mathbb{Z}_3^{2^k}.$
\end{conjecture}
\begin{conjecture}\label{Conjecture 8.4}
$w^1_{\A_3}(G_{t,s^k})=k$
\end{conjecture}

\begin{tabular}{|l| c| c| c| c| c|}\hline
 k & $G_{t,s^k}$ & $torH_{\A_3} ^{1,2v-6}$& $torH_{\A_3} ^{1,2v-5}$& $torH_{\A_3} ^{1,2v-4}$& $H_{\A_3} ^{1,2v-3}$
 \\ \hline

 1&\psfig{figure=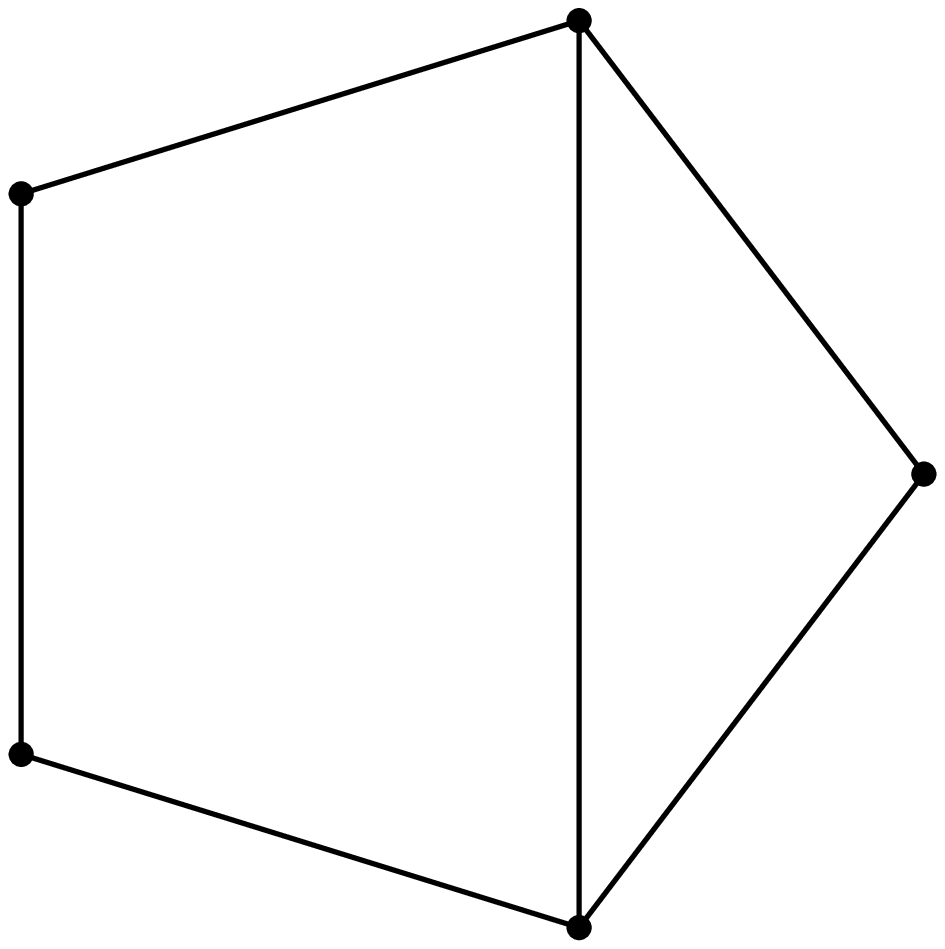,height=1.5cm} &0&0&$\mathbb{Z}_3^2$&$\mathbb{Z}_3 \oplus \mathbb{Z}$  \\ \hline
 2&\psfig{figure=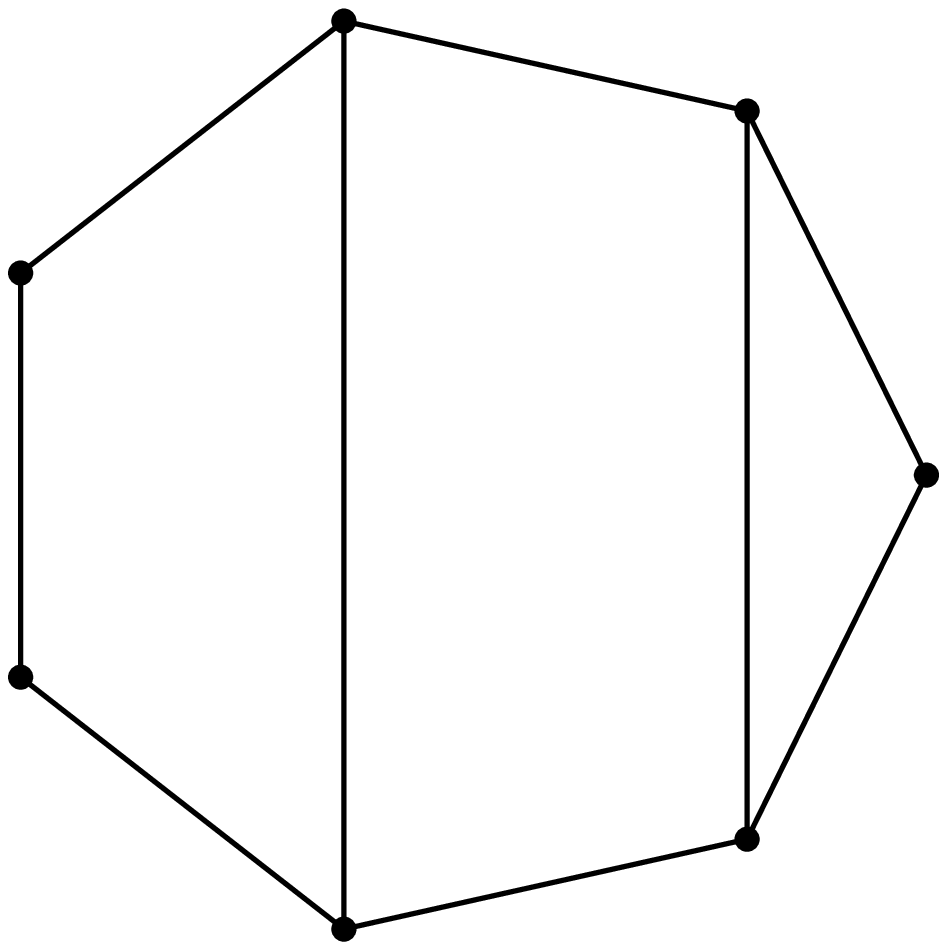,height=1.5cm} &0&$\mathbb{Z}_3^4$&$\mathbb{Z}_3^4$&$\mathbb{Z}_3 \oplus \mathbb{Z}^2$ \\ \hline
 3&\psfig{figure=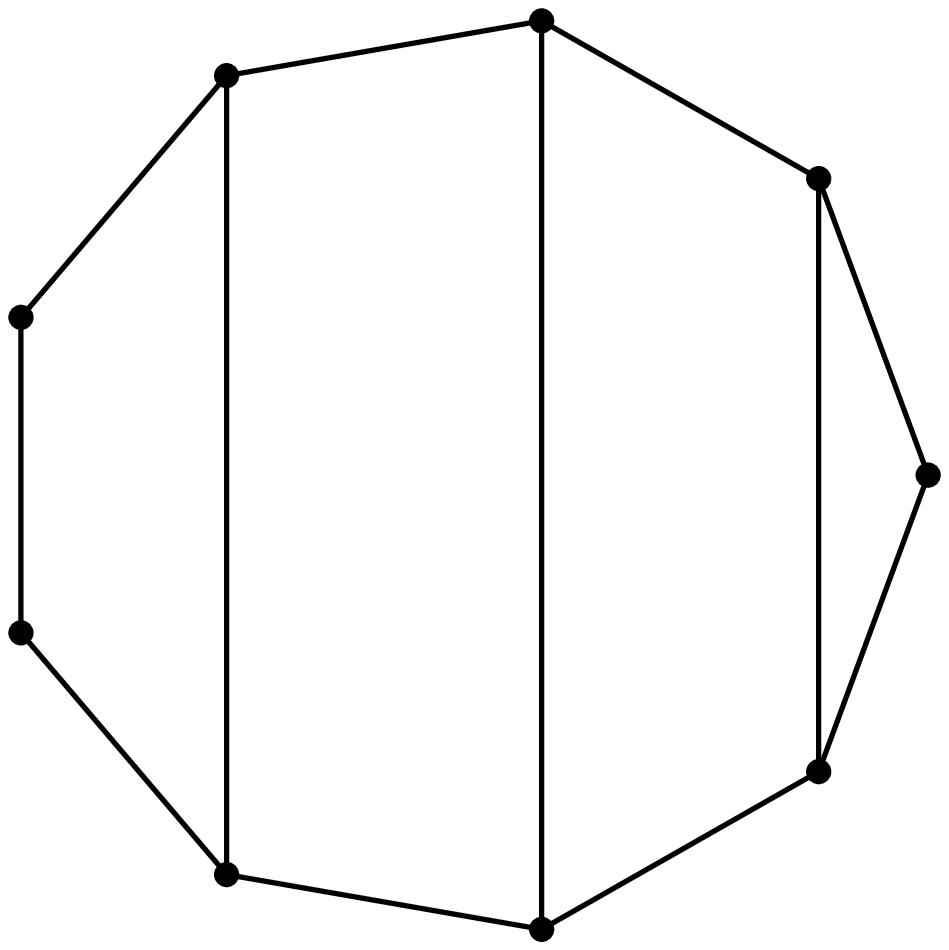,height=1.5cm} &
$\mathbb{Z}_3^{8}$&$\mathbb{Z}_3^{12}$&$\mathbb{Z}_3^{6}$&
$\mathbb{Z}_3 \oplus \mathbb{Z}^3$ \\ \hline
 4&\psfig{figure=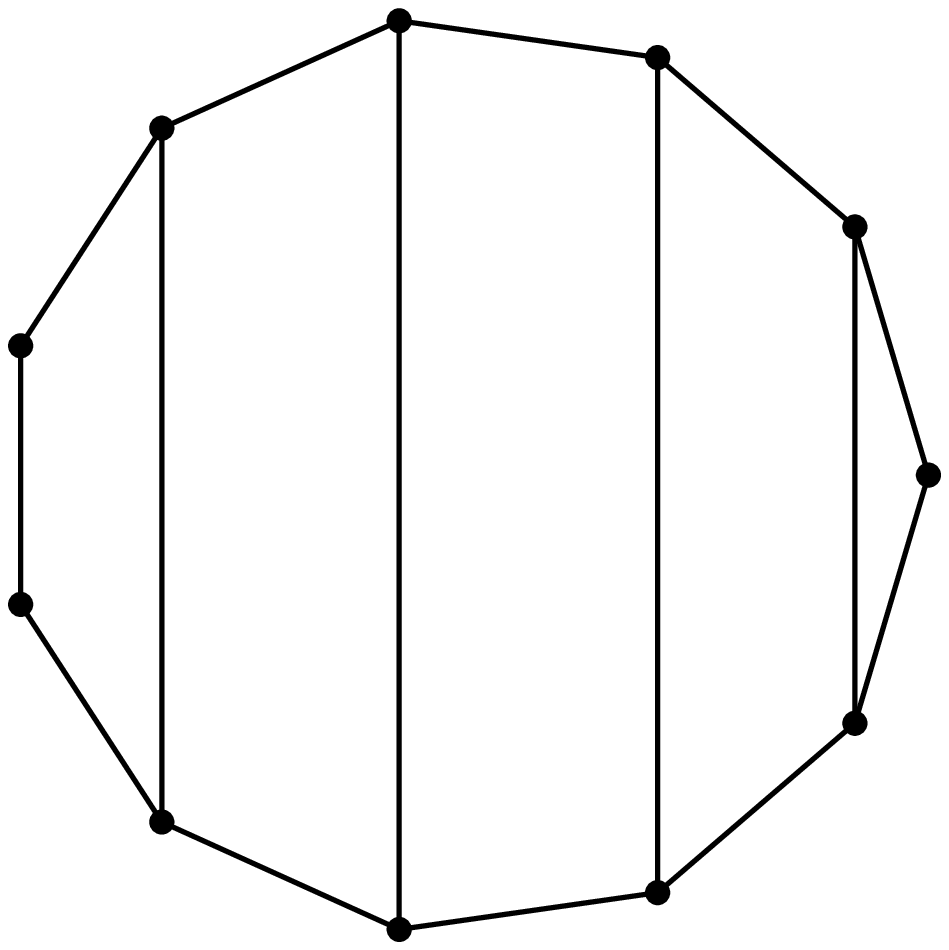,height=1.5cm} & 
$\mathbb{Z}_3^{32}$&$\mathbb{Z}_3^{24}$&$\mathbb{Z}_3^{8}$&$\mathbb{Z}_3 \oplus \mathbb{Z}^4$ \\ \hline
\end{tabular}
\ \\

Next, we consider family of wheels 
$W_k$, that is, cones over $(k-1)$-gones.
$H^{1,2k-3}$ was computed in Corollary 4.11 and computational results
presented in the next Table imply the following conjecture:
\begin{conjecture}\label{Conjecture 8.5}
$w^1_{\A_3}(W_n)=[\frac{n-3}{2}]$
\end{conjecture}

\begin{tabular}{|l| c| c| c| c|}\hline
n & $tor H_{\A_3} ^{1,2v-6}(W_n)$& $torH_{\A_3} ^{1,2v-5}(W_n)$&
$H_{\A_3} ^{1,2v-4}$& $H_{\A_3} ^{1,2v-3}(W_n)$
 \\ \hline

4 &0&0&0&$\mathbb{Z}_2 \oplus\mathbb{Z}_3^3 \oplus \mathbb{Z}^2$  \\
\hline

 5 &0&0&$\mathbb{Z}_3^3$&$\mathbb{Z}_3^3 \oplus
\mathbb{Z}^5$ \\ \hline

 6&0&0&$\mathbb{Z}_3^{5}$&$\mathbb{Z}_2 
\oplus\mathbb{Z}_3^5 \oplus \mathbb{Z}^5$ \\
\hline

7 &0& $\mathbb{Z}_3^{2}$&$\mathbb{Z}_3^{9}$& $\mathbb{Z}_3^5
\oplus \mathbb{Z}^7$ \\ \hline 8

&0&$\mathbb{Z}_3^{7}$&$\mathbb{Z}_3^{14}$&$\mathbb{Z}_2 \oplus
\mathbb{Z}_3^7 \oplus \mathbb{Z}^7$ \\ \hline 9

&$\mathbb{Z}_3^{5}$&$\mathbb{Z}_3^{20}$&$\mathbb{Z}_3^{21}$&$\mathbb{Z}_3^7
\oplus \mathbb{Z}^9 $ \\ \hline
\end{tabular}
\ \\
\ \\

Examples, given above, may suggest that if a simple graph $G$
has a cycle of odd length $n$ then the group
$H_{\A_3}^{1,2(v-1)- \frac{n-1}{2}}(G)$ has
 torsion. However we computed that the only torsion in $H_{\A_3}^{1,j}(K_4)$
is supported by $j=5$ and the only torsion in $H_{\A_3}^{1,j}(K_5)$
is supported by $j=7$, despite the fact that $K_4$ and $K_5$ have 5-cycles.
In fact, calculations of $tor H^{1,j}_{\A_3}(K_n)$ suggest:

\begin{conjecture}\label{Conjecture 8.6}  
$w^1_{\A_3}(K_n)= 0$  
\end{conjecture}

\subsection{Computation over the $\A_5$ algebra: $H_{\A_5}^{1, 4v-7}(G)$}
So far, according to calculations for different series of graphs (wheels,
broken wheels and complete graphs) the only torsion is $\Z_2$ and $\Z_5$.
 Furthermore, the rank of $\Z_5$ torsion
is always the number of 3-cycles in a graph, while, for complete graphs,
the rank of $\Z_2$ is the number of 4-cycles in a graph. We should stress
that these results are not sufficient for general conjectures- even for $m=3$
 the above series of
 graphs have only $\Z_3$ and $\Z_2$ torsion but as we proved in Section 
5 in $\A_3$- cohomology any torsion is possible.
\ \\ \
First we consider a family of wheels, $W_n$ and broken wheels:
$W_n^{out}$ is obtained from $W_n$ by deleting an edge from the polygonal 
part of the wheel (treated as a cone over ($n-1$-gon), 
$W_n^{in}$ is obtained from $W_n$ by deleting a ``spike" edge, 
Figure 8.3.
\ \\ \

 \begin{conjecture}\label{Conjecture 8.7}
For $n>4$ the following 
holds\footnote{$H^{1,9}_{\A_5}(W_3^{in})=\mathbb{Z}_5^2 \oplus
\mathbb{Z}$,
$H^{1,13}_{\A_5}(W_4^{in})=\mathbb{Z}_5^2 \oplus \mathbb{Z}^2$\\
$H^{1,9}_{\A_5}(W_3)=\mathbb{Z}_5^2 \oplus \mathbb{Z}_2 \oplus
\mathbb{Z}^2$, $H^{1,13}_{\A_5}(W_4)=\mathbb{Z}_5^4 \oplus
\mathbb{Z}^5$.}:
 $$H^{1,4n-3}_{\A_5}(W_n^{out})=\mathbb{Z}_5^{n-1} \oplus
\mathbb{Z}^{n-2}$$
$$H^{1,4n-3}_{\A_5}(W_n)=\mathbb{Z}_5^{n} \oplus
\mathbb{Z}^{n}$$
$$H^{1,4n-3}_{\A_5}(W_n^{in})=\mathbb{Z}_5^{n-2}
\oplus \mathbb{Z}^{n-2}$$
\end{conjecture}
\ \\
\ \\
\centerline{\psfig{figure=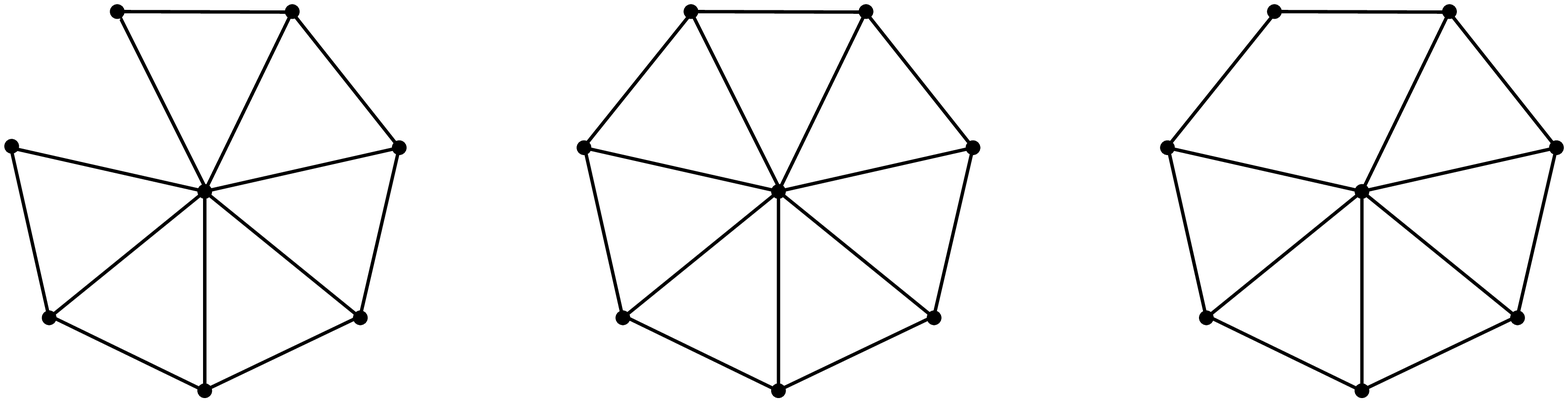,height=2.7cm}}
\begin{tabular}{c c c c l}
\ \ \ $W_8^{out}=W_8- e_{out}$ &\ \ \ \ \ \ \ \ \ \ \ \ \ \ 
&$W_8$& \ \ \ \ \ \ \ \ \ \ \ \ \ \ \ \  &$W_8^{in}=W_8-e_{inside}$ \\
\end{tabular}
\ \\

\centerline{ Figure 8.3:  Wheels and broken wheels } \ \\
\begin{conjecture}\label{Conjecture 8.8}
For algebra $\A_5$ and complete graphs $K_n$ with $n$ vertices the
following is true:
$$ H_{\A_5}^{1,4n-7}(K_n)=\mathbb{Z}_5^{n\choose 3 }
\oplus \mathbb{Z}_2^{n\choose 4} \oplus
\mathbb{Z}^{2 {n\choose 4}}.$$
\end{conjecture}

\subsection{Calculation of $H_{\A_m}^{1,(m-1)(v-2)+1}(G)$}\ \\
Based on computational results we conjecture that the cohomology
$H_{\A_m}^{1,,3m-2}(G)$
for several graphs including a square with a diagonal 
($P_3|P_3$), a ``house" ($P_3|P_4$),  $K_4$
and $K_5$.
\begin{conjecture}\label{Conjecture 8.9}
 For $m >2$ we 
have\footnote{$H^{1,3}_{\A_2}(P_3|P_3)= \Z_2\oplus \Z$.}:
$H^{1,2m-1}_{\A_m}(P_3|P_3)= \Z_m^2 \oplus \Z $.
\end{conjecture}
We checked the conjecture for $m\leq 15$

\begin{conjecture}\label{Conjecture 8.10}
 For $m \geq 2$ we
have $H^{1,3m-2}_{\A_m}(P_3|P_4)= \Z_m \oplus \Z $.
\end{conjecture} 
We checked the conjecture for $m\leq 14$.

\begin{conjecture}\label{Conjecture 8.11}
$(\forall m \geq 4) \  H^{1,2m-1}_{\A_m}(K_4)=\mathbb{Z}_m^4
\oplus \mathbb{Z}_2 \oplus \Z^2.$
\end{conjecture}
We checked the conjecture for $m\leq 14$.
\begin{conjecture}\label{Conjecture 8.12}
$H^{1,3m-2}_{\A_m}(K_5)=\left\{%
\begin{array}{ll}
   \mathbb{Z}_2 \oplus  \mathbb{Z}^{5}, & \hbox{m=2;} \\
    \mathbb{Z}_2 \oplus  \mathbb{Z}_3^{4}
\oplus \mathbb{Z}^{10}, & \hbox{m=3;} \\
\mathbb{Z}_2^{5} \oplus \mathbb{Z}_m^{10} \oplus \mathbb{Z}^{10},
& \hbox{m=4;} \\
    \mathbb{Z}_2^{5} \oplus \mathbb{Z}_m^{10} \oplus \mathbb{Z}^{10}, & m>4, m \ odd; \\
    \mathbb{Z}_2^{10} \oplus
 \mathbb{Z}_m^{10} \oplus  \mathbb{Z}^{10}, & m>5,m \ even. \\
\end{array}%
\right.$
\end{conjecture}
\ \\
We checked the conjecture for $m\leq 9$.

We hope that conjectures and calculations listed in this Section will encourage future research.
\ \\
\ \\

\ \\ \ \\ \ \\
\noindent \textsc{Dept. of Mathematics, Old Main Bldg., 1922 F St.
NW
The George Washington University, Washington, DC 20052}\\
e-mails: {\tt pabiniak@gwu.edu}, {\tt przytyck@gwu.edu},
         {\tt radmila@gwu.edu}
\end{document}